\documentclass[a4paper,12pt]{amsart}\usepackage{amsfonts,amssymb,bbm}
\usepackage{enumerate,mathrsfs,xfrac,xspace,hyperref,rotating} \UseRawInputEncoding
\def\C{\Bbb{C}}
\def\k{\Bbbk}
\def\K{\mathbbm{K}}

\def\N{\Bbb{N}}\def\Q{\Bbb{Q}}\def\R{\Bbb{R}}

\def\bl{\langle}\def\br{\rangle}

\def\capl{\mathop\cap\limits}

\newcommand{\quot}[2]{{\footnotesize\left.\raisebox{1.2ex}{$#1$}\!\! \ensuremath\diagup \!\!\raisebox{-1.2ex}{$#2$}\right.}}
\newcommand{\quots}[2]{{\footnotesize\left.\raisebox{0.4ex}{$#1$}\! / \!\raisebox{-0.4ex}{$#2$}\right.}}

\renewcommand{\stackrel}[2]{\ \lower 0.34ex \hbox{$\mathrel{\mathop{#2}\limits^{#1}}$}\ }

\def\tc{\tilde{c}}\def\tf{{\tilde{f}}}\def\tF{\tilde{F}}

\def\tk{{\tilde{\k}}} 

\def\tR{{\tilde{R}}}
\def\tt{{\tilde{t}}}

\def\tx{{\tilde{x}}}\def\tX{{\tilde{X}}}\def\tY{{\tilde{Y}}}
 \def\tz{{\tilde{z}}}\def\tZ{{\tilde{Z}}}

\def\hC{\hat{C}}
\def\hc{\hat{c}}
\def\hf{\hat{f}}\def\hg{{\hat{g}}}
\def\hk{{\hat\k}}\def\hR{{\widehat{R}}}

\def\hX{\hat{X}}\def\hy{\hat{y}}\def\hY{\hat{Y}}\def\hz{{\hat{z}}}
\def\hPhi{\hat{\Phi}}\def\hX{{\hat{X}}}

\def\Ga{\Gamma}\def\tGa{\tilde\Ga}\def\be{\beta}
\def\ep{\epsilon}

\def\cA{\mathscr A}\def\ca{\mathfrak a}\def\tca{\tilde{\ca}}
\def\cC{\mathscr C}

\def\cG{\mathscr G}
\def\cK{{\mathscr K}\!}\def\cL{\mathscr L}\def\cR{\mathscr{R}}

\def\cm{{\frak m}}

\newcommand{\ber}{\begin{array}{l}}\newcommand{\eer}{\end{array}}
\newcommand{\bpm}{\begin{pmatrix}}\newcommand{\epm}{\end{pmatrix}}
\newcommand{\bbm}{\begin{bmatrix}}\newcommand{\ebm}{\end{bmatrix}}
\newcommand{\bM}{\begin{matrix}}\newcommand{\eM}{\end{matrix}}
\newcommand{\bee}{\begin{enumerate}}\newcommand{\eee}{\end{enumerate}}
\newcommand{\bei}{\begin{itemize}}\newcommand{\eei}{\end{itemize}}

\def\wrt{with respect to }
\def\sset{\!\subset\!}\def\sseteq{\!\subseteq\!}\def\ssetneq{\!\subsetneq\!}
\def\spset{\!\supset\!}\def\spseteq{\!\supseteq\!}
\def\smin{\!\setminus\!}

\def\Maps{{\rm Maps}\!\left(X,Y\right)}\def\MapX{{\rm Maps}\!\left(X,(\k^p,o)\right)}
\def\Mapk{{\rm Maps}\!\left((\k^n,o),(\k^m,o)\right)} \def\MapR{{\rm Maps}\!\left((\R^n,o),(\R^m,o)\right)}
\def\MapC{{\rm Maps}\!\left((\C^n,o),(\C^m,o)\right)}

\def\APG{$\mathop \text{AP\!.}\!\cG$\xspace}\def\APR{$\mathop \text{AP.}\cR$\xspace}\def\APL{$\mathop \text{AP.}\cL$\xspace}\def\APLR{$\mathop \text{AP.}\cL\cR$\xspace}
\def\APGam{$\mathop \text{AP.}\Ga$\xspace}
\def\APPGam{$\mathop \text{APP\!.}\!\Ga$\xspace}

\def\APGamP{$\mathop \text{APP\!.}\!\Ga$\xspace}\def\APPG{$\mathop \text{APP\!.}\!\cG$\xspace} \def\APPR{$\mathop \text{APP\!.}\!\cR$\xspace}
\def\APPL{$\mathop \text{APP\!.}\!\cL$\xspace}\def\APPC{$\mathop \text{APP\!.}\!\cC$\xspace}\def\APPK{$\mathop \text{APP\!.}\!\cK$\xspace}
\def\APPLR{$\mathop \text{APP\!.}\!\cL\cR$\xspace}

\def\SAPG{$\mathop \text{SAP\!.}\!\cG$\xspace}\def\SAPGam{$\mathop \text{SAP\!.}\!\Ga$\xspace} \def\SAPR{$\mathop \text{SAP\!.}\!\cR$\xspace} \def\SAPL{$\mathop \text{SAP\!.}\!\cL$\xspace}  \def\SAPC{$\mathop \text{SAP\!.}\!\cC$\xspace} \def\SAPK{$\mathop \text{SAP\!.}\!\cK$\xspace}\def\SAPLR{$\mathop \text{SAP\!.}\cL\cR$\xspace}

\def\RmX{R^{\oplus m}_X}\def\RmXY{R^{\oplus m}_{X\!\times\! Y}}

\newcommand{\beq}{\vspace{-0.08cm}\begin{equation}}\newcommand{\eeq}{\vspace{-0.08cm}\end{equation}}  %{
\newtheorem{Lemma}{Lemma}[section]\newcommand{\bel}{\vspace{-0.1cm}\begin{Lemma}}\newcommand{\eel}{\vspace{-0.1cm}\end{Lemma}}
\newtheorem{Example}[Lemma]{Example}\newcommand{\bex}{\vspace{-0.1cm}\begin{Example}\rm}%\newcommand{\eex}{\hfill$\lrcorner$\end{Example}}
\newcommand{\eex}{\vspace{-0.1cm}\end{Example}}
\newtheorem{Proposition}[Lemma]{Proposition}\newcommand{\bprop}{\begin{Proposition}}\newcommand{\eprop}{\end{Proposition}}
\newtheorem{Property}[Lemma]{Property}\newcommand{\bproperty}{\begin{Property}}\newcommand{\eproperty}{\end{Property}}
\newtheorem{Definition-Proposition}[Lemma]{Definition-Proposition}

\def\bpr{~\\{\em Proof.\ }}
\newcommand{\epr}{{\hfill\ensuremath\blacksquare}\\}
\newtheorem{Theorem}[Lemma]{Theorem}\newcommand{\bthe}{\vspace{-0.1cm}\begin{Theorem}}\newcommand{\ethe}{\vspace{-0.1cm}\end{Theorem}}
\newtheorem{Definition}[Lemma]{Definition}\newcommand{\bed}{\vspace{-0.1cm}\begin{Definition}}\newcommand{\eed}{\end{Definition}\vspace{-0.1cm}}
\newtheorem{Remark}[Lemma]{Remark}\newcommand{\beR}{\vspace{-0.1cm}\begin{Remark}\rm}\newcommand{\eeR}{\vspace{-0.1cm}\end{Remark}}
\newtheorem{Corollary}[Lemma]{Corollary}\newcommand{\bcor}{\vspace{-0.1cm}\begin{Corollary}}\newcommand{\ecor}{\vspace{-0.1cm}\end{Corollary}}

\newcommand{\bet}{\begin{tabular}{lll}}\newcommand{\eet}{\end{tabular}}

\title[]{F\MakeLowercase{urther results on} A\MakeLowercase{rtin approximation,}\\ \MakeLowercase{for group-actions on mapping-germs}
  $\Maps$\ \MakeLowercase{and for quivers of maps}}
\author[]{D\MakeLowercase{mitry} K\MakeLowercase{erner}}
\address{Department of Mathematics, Ben Gurion University of the Negev, P.O.B. 653, Be'er Sheva 84105, Israel. dmitry.kerner@gmail.com}
\date{\today\ \  filename: \jobname.tex}
\thanks{I was supported by the Israel Science Foundation, grants No.  1910/18 and 1405/22}


\begin{document} 
\maketitle
\vspace{-1cm}
\begin{abstract}
Consider (analytic, resp. algebraic) map-germs, $\Mapk.$ These germs are traditionally studied up to the right ($\cR$), let-right ($\cL\cR$) and contact ($\cK$) equivalences.  Below $\cG$ is one of these groups. 

 An important tool in this study is the Artin approximation: any formal equivalence, $\tf\stackrel{\hat\cG}{\sim}f,$ is approximated by ordinary (i.e. analytic, resp. algebraic) equivalence,
$\tf\stackrel{\cG}{\sim}f.$

We consider maps of (analytic, resp. algebraic) scheme-germs, with arbitrary singularities, $\Maps,$ and establish stronger versions of this property (for $\cG$): the Strong Artin approximation and the P\l oski approximation.

As a preliminary step we study the contact equivalence for maps with singular targets.

\medskip

In many cases one works with multi-germs of spaces, and with their ``muti-maps". More generally, ``quivers of map-germs" occur in various applications. The needed tools are the Strong Artin approximation for quivers and the P\l oski version. We establish these for directed rooted trees.
\end{abstract}

\setcounter{secnumdepth}{6} \setcounter{tocdepth}{1}\vspace{-0.7cm}
 \tableofcontents

\vspace{-1cm}
\section{Introduction}
The classical Artin approximation, \cite{Artin.68},\cite{Artin.69}, reads: given a system of algebraic (resp. analytic) equations of implicit function type,
 $F(x,y)=0,$ every formal solution, $F(x,\hy(x))=0,$  is approximated by algebraic (resp. analytic) solutions,
 $F(x,y(x))=0$.  See \S\ref{Sec.Background.AP} for more detail.

 The subject ``Artin approximation" has been studied extensively, see e.g. \cite{Popescu}, \cite{Rond} for numerous results and further references. It became an everyday tool in Singularity Theory, Local Algebraic/Analytic Geometry, and in Commutative Algebra.

 The known versions of Artin-approximation are very general, while for particular applications one needs more specialised statements.
 In some cases such ``specializations" are not immediate or just fail.

\subsection{Approximation problems for mapping-germs}
A central object of study in Singularity Theory (and numerous related areas) are the mapping germs $\Mapk.$
 Here
 \bei
 \item $\k$ is a field of arbitrary characteristic. In the analytic case $\k$ is a complete normed field.
 \item $(\k^n,o),(\k^m,o)$ are the (algebraic/analytic/formal) germs of affine spaces.
  \eei
More generally, one considers  $\Maps,$ where $X\sseteq(\k^n,o), $  $Y\sseteq(\k^m,o) $ are  the (algebraic/analytic/formal) subgerms, with arbitrary singularities. Various questions of deformation theory force $\k$ to be a (local) ring, rather than a field.

 \medskip

These maps are classically studied up to the  right ($\cR$), left ($\cL$), left-right ($\cL\cR$) and contact ($\cK$) group-actions.
 E.g. the group action $\cL\cR\circlearrowright\Maps$ goes by $f\rightsquigarrow \Phi_Y\circ f\circ\Phi_X^{-1},$ where
  $\Phi_X\in Aut_X$ and $\Phi_Y\in Aut_Y$ are automorphisms of the source and of the target.
 See \S\ref{Sec.Background.R.LR.equiv} and \S\ref{Sec.Background.Contact.Equivalence} for more detail.

Below $\cG$ is one of these groups $\cR,\cL,\cL\cR,\cK.$ The first fundamental questions are:
\bee[\!\!$\bullet$]
\item
To decide whether two (algebraic/analytic/formal) maps, $X\!\stackrel{f,\tf}{\to}\!Y,$ are $\cG$-equivalent, i.e. $\tf\!\in\! \cG f.$
\item
To bring a map $X\stackrel{f}{\to}Y $  to a particular ``normal" form. E.g. to eliminate monomials of certain type in the power series $f(x),$ or to find a $\cG$-representative of $f$ inside a particular subspace/substratum of $\Maps.$
\item To reduce the study of deformations of the map $X\stackrel{f}{\to}Y $ to the ``tangent level" (i.e. deformations over $\quots{\k[\ep]}{(\ep)^2}$) or to the level of obstruction functors ($T^1,T^2,\dots$)
\eee
In ``simple" cases this can be done by geometric arguments, or by studying the lowest order monomials of $f.$ In more complicated cases one can only invoke inductive arguments, eliminating monomials order-by-order. Such arguments give formal solutions  (and even this can be non-obvious for the left-right equivalence $\cL\cR$). 
 But one needs analytic/algebraic solutions. Hence the need in the corresponding Artin approximation.
\bee[\!\!$\bullet$]
\item
For $\cR$-equivalence the condition to resolve is an  implicit function equation, $\tf(x)=f(\phi(x))$ with the unknown $\phi.$ Here the standard Artin approximation readily applies.

\item For $\cL\cR$-equivalence (and for $\cK$-equivalence, when the target $Y$ is non-smooth) the involved equations are not of implicit-function-type. Moreover, the relevant Artin approximation, \APLR, does not hold in the analytic case,  \cite[Fact.1.4]{Shiota.1998}, \cite[pg. 1061]{Shiota.2010}.
\eee

\medskip

    The first version of \APLR was established in \cite{Shiota.1998} for $\MapR$  that are either Nash or [analytic and of ``finite singularity type"].
 In \cite{Kerner.LRAP} the approximation \APLR was extended to an arbitrary base-field $\k$, for $\Maps$  that are either $\k$-Nash  or [$\k$-analytic and ``of weakly finite singularity type"].

\medskip

For the further study of $\Maps$, e.g. in \cite{Kerner.Unfoldings} and \cite{Kerner.Change.of.Base.Field}, we need stronger versions of Artin approximations. This is the goal of the current paper. Though it is the natural continuation of \cite{Kerner.LRAP}, the results and proofs are independent.

\subsection{The results} Below $\cG$ is one of the groups $\cR,\cL,\cL\cR,\cC,\cK,$ see \S\ref{Sec.Background.R.LR.equiv} and  \S\ref{Sec.Background.Contact.Equivalence}.

 \subsubsection{Strong Artin approximation for $\cG$-equivalence} Occasionally (when maps $X\stackrel{f,\tf}{\to}Y$ are complicated) one cannot even establish a formal equivalence, $\tf\stackrel{\hat\cG}{\sim}f.$
  The only possible result is an approximate solution, $\tf\stackrel{mod\ \cm^d}{\equiv}g_d f$ for some $d\gg1$ and $g_d\in\cG.$
  To establish the ordinary (e.g. analytic/algebraic) equivalence, one needs the ``Strong Artin $\cG$-approximation":
 \\{\bf Theorem \ref{Thm.SAP.G.and.APP.G}, {\pmb\SAPG}.}  (roughly) \ Given two maps $X\stackrel{f,\tf}{\to}Y,$ there exists a function $\N\stackrel{\be}{\to}\N$ ensuring:
   if  $\tf \equiv g_d f\ mod\ \cm^{\be(d)}$ for some $g_d\in \cG,$ then $\tf=g f$ for a corresponding $g\in \cG$ that satisfies
     $g \equiv g_d\ mod\ \cm^d .$

This can be restated topologically. The set of all maps, $\Maps,$ has the natural $\cm$-adic filtration. And \SAPG insures: the orbit $\cG f\sset \Maps$ is closed in this filtration-topology. See \S\ref{Sec.GSAP.GAPP.definitions} for detail.

\subsubsection{P\l oski version of Artin approximation}
 The classical theorem of P\l oski (extended by Popescu and others, \S\ref{Sec.Background.APP}) reads: given a formal solution, $F(x,\hy(x))=0,$ there exists a parameterized (analytic/algebraic) solution, $F(x,y_z(x))=0,$ such that $\hy(x)=y_{\hz(x)}(x)$ for a certain specialization $\hz(x).$
  Thus the formal solutions are ``artificial" and ``not inherent". We give
  the corresponding P\l oski-type property (\APPG) for $\cG$-equivalence of maps.
\\{\bf Theorem \ref{Thm.SAP.G.and.APP.G}, {\pmb\APPG}.}  (roughly) \ Given a map $X\stackrel{f}{\to}Y,$ and a formal equivalence
     $\tf=\hg f,$ there exists a parameterized equivalence $\tf=g_z f$ such that $\hg=g_\hz$ for certain specialization $\hz.$

\subsubsection{Remarks on Theorem \ref{Thm.SAP.G.and.APP.G}}\label{Rem.Introduction}
\bee[\hspace{-0.3cm}\bf i.]

\item   While the classical group-actions $\cG\circlearrowright\MapX,$ for $\cG=\cR,\cL,\cL\cR,\cK,$ are well-studied, we could not find references for the actions
 $\cG\circlearrowright\Maps$ with singular target. The relevant basic notions are set-up in \S\ref{Sec.Background.Contact.Equivalence}, and studied further in \cite{Kerner.Change.of.Base.Field}.

\item
For the group $\cG=\cR$ the approximations \SAPR,\APPR are immediate example of the classical properties SAP, APP.
 Indeed, in this case the condition $\tf\stackrel{\cR}{\sim}f$ is an equation of implicit function type.

  But for $\cG=\cL\cR$ and for [$\cG=\cK,$ with singular target, $Y\ssetneq (\k^m,o)$] the involved conditions are not implicit-function-equations. The proof involves a little trick, used e.g. in \cite{Shiota.1998}, transforming the condition $\tf\stackrel{\cG}{\sim}f$ into an implicit-function-equation with nested solutions.

\item The properties \SAPG and \APPG approximate the solutions $g_d,\hg$ by ordinary solutions \wrt the filtrations of local rings by their maximal ideals, $\cm^\bullet_X\sset R_X,$
    $\cm^\bullet_Y\sset R_Y.$ Geometrically, the approximation holds ``near the points $V(\cm_X)\in X$, $V(\cm_Y)\in Y$."

    For various applications one needs stronger approximations, ``along subloci of $X,Y$". E.g. given a formal solution $\tf=\hg f,$ we look for an ordinary  solution $\tf=g_j f$ that satisfies: $\hg \equiv g_j\ mod\ I^j.$ (Here $\sqrt{I}\ssetneq\cm_X,$ i.e. the locus $V(I)\sset X$ is of positive dimension.) This is impossible, as even the classical Strong Artin Approximation does not hold for the filtration $I^\bullet.$ Yet we approximate the solutions $\hg,g_d$ in a related filtration: $\hg \equiv g_j\ mod\ \ca\cdot \cm^j.$ See \S\ref{Sec.SAP.G.APP.G.remarks} for this version of \SAPG and \APG.

\item  In various applications one uses the ``filtration-unipotent" subgroups $\cG^{(j)}<\cG,$ see e.g. \cite{B.K.motor}, \cite{BGK.20}, \cite{Kerner.Group.Orbits}. Roughly speaking, the elements of $\cG^{(j)}$ are identities modulo ``higher order terms". 
   Hence the need in the relevant approximation.  In \S\ref{Sec.SAP.G.APP.G.unipotent} we deduce the properties SAP and APP for  $\cG^{(j)}.$
\item
This theorem is for  $\k$ a ring (not just a field), which is needed in the study of  unfoldings and deformations.

This and other results of the paper hold for a larger class of rings, in many cases just being local, excellent, henselian suffices.
 But we work just with the rings of formal/analytic/agebraic power series, as these are the most important in applications.
\eee 

\subsubsection{Artin approximations for quivers of maps}\label{Sec.Intro.AP.for.quivers}
 The versions  of Artin approximation for the right, left, left-right equivalences  (\APR,      \APL, \APLR) can be encoded by the simplest ``quivers of maps":
\beq\label{Eq.quiver.simplest}
\Phi_X \circlearrowright X\stackrel{f}{\to}Y,\quad\quad\quad\quad\quad
 X\stackrel{f}{\to}Y\circlearrowright \Phi_Y,\quad\quad\quad\quad\quad
\Phi_X \circlearrowright X\stackrel{f}{\to}Y\circlearrowleft \Phi_Y.
\eeq
On many occasions (in the study of map-germs) one has to work with multi-germs of maps, $\amalg X_j\stackrel{\{f_j\}}{\to}Y,$ and dually
 $X\stackrel{\{f_j\}}{\to}\amalg Y_j .$ The corresponding $\cL\cR$-approximation problems are encoded by the quivers:
\beq\label{Eq.quiver.multigerms}
\bM \Phi_1\circlearrowright X_1  \\\dots\\ \Phi_k\circlearrowright X_k\eM
\bM \begin{turn}{-30}\stackrel{f_1}{\to}\end{turn} \\ \begin{turn}{30}\stackrel{f_k}{\to}\end{turn}\eM
  Y\circlearrowleft \Phi_Y,\hspace{2cm}
 \Phi_X \circlearrowright X
\bM\begin{turn}{30}\stackrel{f_1}{\to}\end{turn}\\\begin{turn}{-30}\stackrel{f_k}{\to}\end{turn}\eM
\bM Y_1 \circlearrowleft \Phi_1\\\dots\\ Y_k \circlearrowleft \Phi_k\eM\quad .
\eeq
Various other quivers of maps arise naturally in Singularity Theory, Dynamical Systems, Control Theory. E.g. the classical study of self-maps $f\circlearrowright(\k^n,o)$ up to conjugation, $f\sim \Phi\circ f\circ \Phi^{-1},$ is encoded by the quiver
 $\Phi\circlearrowright (\k^n,o)\circlearrowleft f.$ Recall that in this case one has no formal-to-analytic approximation, by the classical results of Poincar\`e, Siegel and others. (See Remark \ref{Rem.Gamma.AP}.)

 All this calls for the general treatment of ``Artin approximation for quivers of map-germs", and its Strong Artin-P\l oski-Popescu versions.
A remark, the classical study of quivers was for {\em linear} maps, \cite{Kirillov}, thus no questions of Artin-approximation were involved.

\medskip

Take a directed graph $\Ga,$ whose vertices, $\{X_v\},$ are scheme-germs, while the edges
 are germs of maps, $X_v\stackrel{f_{wv}}{\to}X_w.$ A morphism of such quivers,
 \beq
  (\Ga,\{\tX_\bullet\}, \{\tX_v\stackrel{\tf_{wv}}{\to}\tX_w\}_{wv})\quad \quad \to\quad \quad
  (\Ga,\{X_\bullet\} , \{X_v\stackrel{f_{wv}}{\to}X_w\}_{wv}),
  \eeq
 is a collection of morphisms of schemes $\{\Phi_v: \tX_v\to X_v\}_v$ that ensure commutativity of diagrams
  ``$\tGa$ over $\Ga$", see \S\ref{Sec.Quivers.Definition}.
  This gives the system of equations $\{\Phi_w\circ \tf_{wv}= f_{wv}\circ \Phi_v\}_{wv}.$
  The maps $\{f_v\},\{\tf_v\}$ are prescribed, while   $\{\Phi_v\}$ are unknowns, in different sets of variables.
  One can impose various additional conditions on $\{\Phi_v\}$, e.g. $\Phi_v$ is invertible,   $\Phi_v=Id_{X_v},$ or $\Phi_v$ is identity modulo higher order terms, or $\Phi_v$ sends a subgerm $\tZ\sset \tX$ to a subgerm $Z\sset X.$

  Theorem \ref{Thm.Gamma-AP.for.k<x>} establishes the ``Strong Artin approximation for $\Ga$" and the ``P\l oski approximation for $\Ga$"  in the algebraic/formal cases, when $\Ga$ is a rooted tree.
 No analytic version of this result is given, as the proof inevitably uses the nested Artin approximation. (Which is absent in the analytic case.)

\subsubsection{Artin approximations for (quivers of) maps with base change}\label{Sec.Intro.AP.for.quivers.parameters}
The quivers considered in \S\ref{Sec.Intro.AP.for.quivers} are over a fixed base $Spec(\k),$ and the maps correspond to morphisms of $\k$-algebras,
  To study induced deformations, unfoldings and normal forms (of maps and of other objects) it is important to compare objects over various bases ($Spec(\tk)$ vs $Spec(\k)$), and to consider morphisms of rings involving homomorphisms  $\k\to \tk.$

 Thus  we consider families of scheme-germs, $X_t$, and  of maps-germs $X_t\stackrel{f_t}{\to} Y_t.$
  Accordingly we have the property ``\APG  with change of parameters".

 More generally, we have quivers of families of maps,
  $ (\Ga,\{X_{t,v}\}_v, \{X_{t,v}\!\stackrel{f_{t,wv}}{\to}\!X_{t,w}\}_{wv}).$ Corollary \ref{Thm.Gamma-AP.with.parameters}
   gives the relevant versions  of approximation: ``SAP.$(\Ga,t)$" and APP.$(\Ga,t)$.
   
   An immediate application is Example \ref{Ex.normal.form.unfolding}, bringing unfoldings to normal forms. 
%\vspace{-0.2cm}

\subsection{Further applications} Theorem \ref{Thm.SAP.G.and.APP.G} is used in \cite{Kerner.Change.of.Base.Field} to trace the $\cG$-equivalence of maps under the extension of the base field, $\k\sset \K.$ We get (among other results) an unexpected property: If $\tf\in \cG^{(1)}_\K f$ then
  $\tf\in \cG^{(1)}_\k f.$

The results of \S\ref{Sec.Intro.AP.for.quivers.parameters} are used in \cite{Kerner.Unfoldings} to extend the classical theory of unfoldings of $\MapR$ and $\MapC$ to $\Maps$ in arbitrary characteristic.

\subsection{Acknowledgements} Thanks to Dorin Popescu and Guillaume Rond for useful advice.

% \subsection{Structure of the paper}

\section{Definitions and  auxiliary results}

\subsection{Notation  and conventions (Germs of schemes and their morphisms)}\label{Sec.Background.Schemes.Maps}
 For the general introduction to  $\Mapk$
 in the real/complex-analytic case see \cite{Mond-Nuno}.
  For the definitions and results over arbitrary fields see \cite{Kerner.Group.Orbits}.

 \medskip 
 
   We use   multivariables, $x=(x_1,\dots,x_n),$  $y=(y_1,\dots,y_m),$ $t=(t_1,\dots,t_r),$ $F=(F_1,\dots,F_r).$
 The notation  $F(x,y)\in (x,y)^N $  means: every component of the vector $F$  satisfies this condition.

\subsubsection{}\label{Sec.Background.Rings} The source and the target of maps are (formal, resp. analytic,   algebraic) scheme-germs:
\bei
\item  $X\!=\!Spec(R_X)$, where  $R_X\!=\!\quots{\k[\![x]\!]}{J_X},$ resp. $R_X\!=\!\quots{\k\{x\}}{J_X},$    $R_X\!=\!\quots{\k\bl x\br}{J_X}.$
\\$Y=Spec(R_Y)$, where $R_Y\!=\!\quots{\k[\![y]\!]}{J_Y},$ resp. $R_Y\!=\!\quots{\k\{y\}}{J_Y},$     $R_Y\!=\!\quots{\k\bl y\br}{J_Y}.$

Occasionally we use the germ $X\times Y.$ Its local ring $R_{X\!\times\! Y}$ is one of $\quots{\k[\![x,y]\!]}{J_X+J_Y},$
 $\quots{\k\bl x,y\br}{J_X+J_Y},$ $\quots{\k\{x,y\}}{J_X+J_Y}.$

\item In the formal case, i.e. $\quots{\k[\![x]\!]}{J_X},$ \ $\k$ is a complete local ring.

 \medskip

In the analytic case (i.e. $\quots{\k\{x\}}{J_X}$): \  $\k=\quots{\k_o\{t\}}{\ca},$ where $\k_o$ is either a normed field  (complete \wrt its norm), or a Henselian DVR over a field.
 (The norm is always non-trivial.) The simplest cases are   $\k=\k_o $ and   $\k=\k_o\{t\} .$

% The \!main relevant examples of \!DVR are  the rings
% $\k_o[\![\tau]\!], \k_o\bl \tau\br, \k_o\{\tau\}.$ (Here  $\k_o$ is a field.)

In the algebraic case (i.e. $\quots{\k\bl x\br}{J_X}$):  \  $\k$ is any  field or an excellent Henselian local ring over a field.
  The simplest case is $\k=\quots{\k_o\bl t\br}{\ca},$ e.g. $\k=\k_o$ or $\k=\k_o\bl t\br.$
 See page 4 of \cite{Denef-Lipshitz} for more detail.
\\ 
Recall that $\R\bl x\br$ is the ring of Nash power series.
 We call the case $\quots{\k\bl x\br}{J_X}$  by  ``$\k$-Nash" or just ``Nash", to avoid any confusion with the algebraic rings like  $\k[x]$ or  $\k[x]_{(x)}.$
\eei

Through the paper we call $\k$ (resp. $\k_o$) {\em the base field}.

\medskip

  In the $\k$-smooth case, i.e. $J_X=0,J_Y=0,$ we denote  $(\k^n,o):=Spec(R_X)$, $(\k^m,o):=Spec(R_Y).$ These are formal/analytic/$\k$-Nash germs,
     depending on the context.

  Otherwise the singularities of $X,Y$ are arbitrary.

Recall that the rings $R=\k[\![x]\!],\k\{x\},\k\bl x\br$ admit substitutions: if $f\in R$ and $g_1,\dots,g_n\in \cm_X\sset R,$ then $f(g_1,\dots,g_n)\in R.$

  \subsubsection{}
A (formal/analytic/Nash) map $f\in \Maps$ is defined by the corresponding homomorphism of $\k$-algebras, $f^\sharp\in Hom_\k(R_Y,R_X).$ All our homomorphisms are local.
\\In the   case   $J_X=0,J_Y=0 $    one gets the maps of $\k$-smooth germs,  $\Mapk.$

 Fix   embeddings $X\sseteq (\k^n,o),$ $Y\sseteq (\k^m,o)$ and some coordinates $(x_1,\dots,x_n)$ on $(\k^n,o),$ and $(y_1,\dots,y_m)$ on  $(\k^m,o).$
 They are sent to the generators of ideals in $R_X,$ $R_Y.$ Abusing notations we write just  $(x)=(x_1,\dots,x_n)\sset R_X$ and $(y)=(y_1,\dots,y_m)\sset R_Y.$
     If $\k$ is a field, then $(x),(y)$ are the maximal ideals. Otherwise 
     denote the maximal ideal of $\k$  by $\cm_\k.$  Then 
     the maximal ideals are $\cm_X:=\cm_\k+(x)\sset R_X$
      and $\cm_Y:=\cm_\k+(y)\sset R_Y.$  When only $R_X$ is involved we write just $\cm\sset R_X.$

\medskip

A map $X\stackrel{f}{\to} Y$ is represented by
  the tuple of power series, $(f_1,\dots,f_m)\in \cm\cdot\RmX.$
 This embeds the space of maps as a subset $\Maps\sseteq \cm\cdot\RmX.$
 Explicitly, for $R_Y=\quots{\k[\![y]\!]}{J_Y},\quots{\k\{y\}}{J_Y},\quots{\k\bl y\br}{J_Y},$ one has:
 \beq
 \Maps=Hom_\k(R_Y,R_X)=\{f\in \cm\cdot \RmX|\ f^\sharp(J_Y)\sseteq J_X\}.
 \eeq
 If the target is $\k$-smooth, i.e. $Y=(\k^m,o),$  i.e. $R_Y=\k[\![y]\!],$   $\k\{y\},$   $\k\bl y\br$,
  then we get the $R_X$-module   $\MapX= \cm\cdot\RmX.$
 For singular targets the subset $\Maps\sset \cm\cdot \RmX$ is neither additive, nor $\k$-multiplicative.

\subsection{The right, left, and left-right equivalence on $\Maps$}\label{Sec.Background.R.LR.equiv}
 The  (formal/analytic/Nash) automorphisms of $X$ over $Spec(\k)$ are defined via their algebraic counterparts, $\k$-linear automorphisms of the local ring, $Aut_X:=Aut_\k(R_X).$  See e.g. \cite[\S2.4]{Kerner.Group.Orbits}.

 Similarly we denote  $Aut_Y:=Aut_\k(R_Y).$
  In the $\k$-smooth case, $\Mapk,$ these automorphisms are the coordinate changes in the source and the target.

These automorphisms act on the space of maps, defining the right and left group actions:
\bei
\item $\cR:=Aut_{X}: =Aut_\k(R_X)\circlearrowright \Maps$, by $(\Phi_X,f)\to f\circ \Phi^{-1}_X$;
\item $\cL:=Aut_{Y}: =Aut_\k(R_Y)\circlearrowright \Maps$, by $ (\Phi_Y,f)\to \Phi_Y\circ f$.
\eei
These actions define the right, left, and left-right equivalence of maps. E.g.,  $\tf\stackrel{\cL\cR}{\sim}f$ means
 $f=\Phi_Y\circ \tf\circ \Phi^{-1}_X,$ i.e.
$\Phi_Y(\tf(x))=f(\Phi_X(x)),$  for $\Phi_X\in Aut_{(\k^n,o)}$ with $\Phi_X(J_X)=J_X,$ and $\Phi_Y\in Aut_{(\k^m,o)}$ with $\Phi_Y(J_Y)=J_Y.$

 \subsection{The contact equivalence of $\Maps$ and its linearized version}\label{Sec.Background.Contact.Equivalence}
 \bed
 The subgroup $\cC< Aut_{X\!\times\! Y}$ consists of automorphisms that preserve the fibres of the projection $X\times Y\to X$
 and preserve ``the zero-level hyperplane" $X\times\{o\}\sset X\times Y.$ 
 \eed
 Therefore $\cC$ acts as $Id_X$ on $X\!\times\!\{o\},$
  and restricts to automorphisms of the central fibre $\{o\}\!\times\! Y.$

  Explicitly, fix some generators $J_Y=\bl q_\bullet\br,$ then the action $\cC\circlearrowright X\times Y$ is by $(x,y)\to (x,C(x,y)),$
   where
   \beq\label{Eq.group.C.def.equation}
   q_\bullet(C(x,y))=0   \text{ and } C(x,y)=\Phi_Y(y)+(\dots), \quad \text{ with }\quad (\dots)\in (x)\cdot (y)   \text{ and }   \Phi_Y\in Aut_Y.
\eeq

The action on maps, $\cC\circlearrowright\Maps,$ is defined by $f(x)\rightsquigarrow C(x,f(x)).$

\bed
The contact group is the semi-direct product $\cK:=\cC\rtimes Aut_X.$
\eed
 (Note that $\cC\lhd\cK$ is a normal subgroup.) It acts on $X\times Y$ by $(x,y)\to (\Phi_X(x),C(x,y)).$
 The action $\cK\circlearrowright\Maps$ is defined via the graph of $f,$ i.e. $f\rightsquigarrow C(x,f\circ \Phi^{-1}_X).$

 \bex
 When the target is $\k$-smooth, $Y=(\k^m,o),$ this is the classical contact action, see e.g. \cite[pg.108]{Mond-Nuno} or \cite[pg. 156]{AGLV}.
  In this case the space of maps  is an $R_X$-module, $\MapX=\cm\cdot \RmX,$ and the right action $\cR\circlearrowright\MapX$ is $\k$-linear. Namely, $(c\cdot f+\tc\cdot \tf)\circ \Phi^{-1}_X= c\cdot f\circ \Phi^{-1}_X +\tc\cdot \tf \circ \Phi^{-1}_X$.
   The actions $\cL\cR,\cK\circlearrowright Maps(X,(\k^m,o))$ are not $\k$-linear, not even additive.
 \eex

  However, the action  $\cK\circlearrowright\Maps$ admits a ``linearization", even when $Y$ is singular.
\bed\ [The linearized contact groupoid of a map $f\in \Maps$]
\\$\cK_{lin,f}\!:=\{(U,\Phi_X)\!\in\! GL_m(R_X) \rtimes Aut_X  |\ U(x)\!\cdot\! f(x)=C(x,f(x))\ for \ some\ C(x,y)\ as\ above\}.$
\eed

\bex
\bee[\bf i.]
\item For a $\k$-smooth target, $Y=(\k^m,o),$ we get the classical object: $\cK_{lin,f}=GL_m(R_X) \rtimes Aut_X .$  (In particular, $\cK_{lin,f}$ does not depend on $f,$ and is a group.)
\bpr Any $U(x)\in GL_m(R_X)$ acts on $Maps(X,(\k^m,o))$ by $f\to U\cdot f$.  Thus the condition ``$U(x)\cdot f(x)=C(x,f(x))$" is empty.
 Finally, $GL_m(R_X)< \cK_{lin,f}$ is a normal subgroup.
 \epr

 The action $GL_m(R_X) \rtimes Aut_X \circlearrowright Maps(X,(\k^m,o)) $   is $\k$-linear. Namely,
 \[
 (U,\Phi_X)(f+\tf)= U\cdot f\circ \Phi_X^{-1}+U\cdot \tf\circ \Phi_X^{-1}.
\]
     \item
      When the target $Y$ is not $\k$-smooth, $\cK_{lin,f}$ is not a group and depends on $f$. Moreover, $\cK_{lin,f}$ does not act on $\Maps.$
  Namely, given $(U,\Phi)\in \cK_{lin,f}$ and a general  $f\neq \tf\in \Maps,$ the image of $U(x)\cdot \tf(\Phi(x))$ does not lie inside $Y.$
\eee
\eex

In any case, one has the orbit $\cK_{lin,f} f\sseteq \Maps.$
 It  ``linearizes" the group-action $\cK\circlearrowright \Maps$ in the following sense.
\bel\label{Thm.Linearization.of.Contact.Action}
Fix a map $f\!\in\! \Maps.$ Then $\tf\!\stackrel{\cK}{\sim}\!f$ iff $\tf\!\in\! \cK_{lin,f} f .$ \quad
 [Thus $\cK f\!=\!\cK_{lin,f} f.$]
\eel
This linearization is classically known for $\k$-smooth targets,   e.g. \cite[pg.110]{Mond-Nuno} for $Maps((\C^n,o),(\C^m,o)) $ and  \cite[pg.123]{B.K.motor} for $Maps(X,(\k^m,o))$.
\\\bpr $\Rrightarrow$ Let $\tf=C(x,f\circ\Phi^{-1}_X).$ Here $C(x,y)\in (y)\sset R_{X\!\times\! Y},$ thus we can expand: $C(x,y)=A(x,y)\cdot y,$ for some $A(x,y)\in Mat_{m\times m}(R_{X\!\times\! Y}).$
 By definition of $\cC$-group: $C|_{o\times Y}\in Aut_Y,$ therefore $A(x,y)\in GL_m(R_{X\!\times\! Y}).$
 Finally, denote $U(x)=A(x,f\circ\Phi^{-1}_X)\in GL_m(R_X) $ to get:
 $\tf=U(x)\cdot (f\circ\Phi^{-1}_X ).$

 The direction $\Lleftarrow$ is verified similarly.
\epr

\subsection{The classical versions of Artin approximation}\label{Sec.Background.AP} %Below we assume $\hy(x)\in \cm\cdot\hk[\![x]\!].$

\subsubsection{The ordinary version for analytic power series, \cite[Theorem 1.1]{Denef-Lipshitz}}\label{Sec.Background.AP.analytic}
 Let $\k$ be a complete normed field or a Henselian DVR.
 Consider equations $F(x,y)=0$ for $F(x,y)\in (x,y)\sset\k\{x,y\}.$
 Then any formal solution, $\hy(x)\in \cm\cdot \hk[\![x]\!],$ is $\cm$-adically approximated by  analytic solutions.

 \subsubsection{The strong Artin approximation (SAP) for formal power series}\label{Sec.Background.SAP}    \cite{Pfister-Popescu} and \cite[Theorem 7.1]{Denef-Lipshitz}

Let $\k$ be a field or a complete DVR, and take the power series
 $F(x,y)\in (x,y)\sset \k[\![x,y]\!].$ There exists a function $\N\stackrel{\be_F}{\to} \N$ satisfying: for every $d\in \N$ any approximate solution, $F(x,y^{(d)}(x))\in \cm^{\be_F(d)},$ is $\cm^\bullet$-approximated by an ordinary solution, i.e. $F(x,y(x))=0$ with $y(x)- y^{(d)}(x)\in \cm^d.$

\medskip

\noindent A remark: SAP does not hold for the filtration $I^\bullet$ when $\sqrt{I}\ssetneq \cm.$  See e.g. \cite[Example.3.20]{Rond}.

  \subsubsection{The nested AP for algebraic power series,  \cite{K.P.P.R.M.}, \cite{Popescu.86}, \cite[Prop.3.5]{Popescu}, \cite[Theorem 5.8]{Rond}}\label{Sec.Background.AP.nested}
Let $\k$ be a field or
 an excellent Henselian local ring.
 Consider equations $F(x,y)=0$ for $F(x,y)\in (x,y)\sset\k\bl x,y\br.$  Suppose these equations have a  formal  solution  $\hy(x)\in \cm\cdot \k[\![x]\!]$
 that is nested in the following sense:
  \beq\label{Eq.nested.solution.AP}
 \hy_1(x)\in \hk[\![x_1,\dots,x_{i_1}]\!],\ \dots,\ \hy_m(x)\in \hk[\![x_1,\dots,x_{i_m}]\!], \quad \text{ for some } 1\le i_1\le \cdots\le i_m\le n.
 \eeq
   This formal solution    is $\cm$-adically  approximated by  nested algebraic solutions,
    $y_1(x)\in \k\bl x_1,\dots,x_{i_1}\br,$ \dots, $y_r(x)\in \k\bl x_1,\dots,x_{i_m}\br.$

\subsubsection{The nested Strong Artin approximation (SAP) for algebraic power series over a field $\k$,
  \cite{BDLvdD79} and \cite[Corollary 5.11]{Rond}}\label{Sec.Background.SAP.nested} Given $F(x,y)\in (x,y)\sset\k\bl x,y\br,$
  there exists a function $\N\stackrel{\be_F}{\to} \N$ ensuring: for every $d\in \N$ any approximate nested solution, $F(x,y^{(d)}(x))\in \cm^{\be_F(d)},$ as in \eqref{Eq.nested.solution.AP}, is $\cm^\bullet$-approximated by   nested algebraic solutions, i.e. $F(x,y(x))=0$ with $y(x)- y^{(d)}(x)\in \cm^d.$

\noindent{\em The nested SAP for formal power series,} i.e. $F(x,y)\in (x,y)\sset\k[\![x,y]\!],$ holds as before, provided $\k$ is an $\aleph_0$-complete field, \cite[Corollary 16]{Popescu-Rond}. The simplest examples of $\aleph_0$-complete fields are finite fields and uncountable algebraically closed fields, \cite[Theorem 5]{Popescu-Rond}.
 Neither $\bar\Q,$ nor $\R$ are $\aleph_0$-complete.

 \subsubsection{Approximations with parametrized solutions (APP), \cite{Ploski1974}, \cite{Popescu.85}, \cite{Popescu.86}, \cite{Ogoma.94}, \cite{Swan.98}}\label{Sec.Background.APP} Consider (analytic, resp. algebraic) equations $F(x,y)=0.$
  Here    $F(x,y)\in (x,y)\sset \k\{x,y\} $ (with $\k$ as in \S\ref{Sec.Background.AP.analytic}), resp. $F(x,y)\in (x,y)\sset\k\bl x,y\br$
   (with $\k$ an excellent local henselian ring).
   Given a formal solution, $F(x,\hy(x))=0,$ $\hy(x)\in \cm\cdot \k[\![x]\!],$  there exist power series $y(x,z)\in \k\{x,z\}^{\oplus m},$ resp. $\k\bl x,z\br^{\oplus m},$
      and (formal) power series $\hz(x)\in\cm\cdot \k[\![x]\!]^{\oplus c},$ satisfying:
     \beq\label{Eq.Ploski.theorem}
     F(x,y(x,z))=0 \quad \text{ and }\quad  \hy(x)=y(x,\hz(x)) .
     \eeq

     \noindent({\em The nested version of APP for $\k\bl x \br$},   \cite[Theorem 11.4]{Spivakovsky.99}, \cite[Theorem 2.1]{Bilski.Parusiński.Rond})
     If the formal solution $\hy(x)$ is nested, as in \eqref{Eq.nested.solution.AP}, then the parameterizations $y(x,z),$ $ \hz(x)$
     in equation \eqref{Eq.Ploski.theorem} can be chosen nested as well.
     Namely:
     \bei
     \item     $
      y_1(x,z)\in \k\bl x_1,\dots,x_{i_1},z_1,\dots,z_{j_1}\br,$ \dots, $
      y_m(x,z)\in \k\bl x_1,\dots,x_{i_m},z_1,\dots,z_{j_m}\br,$ \quad and
      \item
      $\hz_1(x)\dots \hz_{j_1}(x)\in \k[\![x_1,\dots,x_{i_1}]\!],$
       $\hz_{j_1+1}(x)\dots \hz_{j_2}(x)\in \k[\![x_1,\dots,x_{i_2}]\!],$ \dots,
       \eei
       for some $1\le i_1\le \cdots\le i_m\le m$ and $1\le j_1\le \cdots\le j_m\le c.$

\subsubsection{Remarks on applications of these classical versions}\label{Sec.Background.AP.Remarks}
\bee[\hspace{-0.3cm}\bf i.\!]

\item ({\bf Equations over quotient rings}) Let $R_X$ be one of the rings $\quots{\k[\![x]\!]}{J_X},\quots{\k\{x\}}{J_X},\quots{\k\bl x\br}{J_X}.$ Consider the systems of implicit function equations, $F(x,y)=0.$ Here $F\in R_X[\![y]\!],$ resp. $R_X\{ y\}, R_X\bl y\br.$
    When $J_X\neq0$ this system of equations is understood in the following sense.
     Take any representative $\tF$ of $F$ over $\k[\![x,y]\!],$ resp. $\k\{x,y\} ,\k\bl x,y\br.$
      Then $F(x,y)=0$ means: $\tF(x,y)\equiv0$ mod $J_X,$ i.e. $\tF(x,y)\in J_X.$

Similarly, for $F(x,y)\in \quots{\k[\![x,y]\!]}{J_X+J_Y},\quots{\k\{x,y\}}{J_X+J_Y},\quots{\k\bl x,y\br}{J_X+J_Y}, $ the equation
 $F(x,y)=0$ means  $\tF(x,y)\in J_X+J_Y.$

\item
 ({\bf Presenting the condition $\Phi(I)\sseteq J$ as an implicit function equation})

 Denote by $R$ the rings $\quots{\k[\![x]\!]}{\ca},$   $\quots{\k\bl x\br}{\ca},$ $\quots{\k\{x\}}{\ca}.$
  Fix two ideals $I,J\sset R$
  and  a homomorphism $\Phi:R\to R.$ For the fixed generators, $x=(x_1,\dots,x_n)$ in $R,$ we present $\Phi$ by power series,
   $\{\Phi(x_i)=\tx_i(x)\}_i.$  Fix some (finite set of) generators, $\{q^I_i\}$ of $I,$ and $\{q^J_j\}$ of $J.$
    Then the condition $\Phi(I)\sseteq J$ transforms into the set of equations:
   \beq
  [ \quad \Phi(q^I_i(x))=\joinrel= ]\hspace{1cm}q^I_i(\tx(x))=\sum \tz_{ji}\cdot q^J_j(x),\quad \forall\ i.
   \eeq
This is a system of implicit function equations for the unknowns $\{\tx_i\},\{\tz_{ji}\}$ in $R.$

 Hence the Artin approximation is applicable also to the condition $\Phi(I)\sseteq J.$

\item ({\bf The allowed base-ring  $\k$}) In the classical AP-statements (above) the ring $\k$ is rather restricted, e.g. has to be a field or a DVR. In the statements of \S\ref{Sec.SAP.G.APP.G}, \S\ref{Sec.Quivers}, \S\ref{Sec.AP.Gamma.with.parameters} we use these approximation properties for more general rings, of \S\ref{Sec.Background.Schemes.Maps}.ii.
    This is possible due to the following argument.

    Let the base ring $\k$ be one of $\quots{\k_o[\![t]\!]}{\ca},$
      $\quots{\k_o\{t\}}{\ca},$  $\quots{\k_o\bl t\br}{\ca},$ where $k_o$ is a field or a DVR. Accordingly present $R_X$ as
     $\quots{\k_o[\![t,x]\!]}{\ca\!+J_X},$ resp   $\quots{\k_o\{t,x\}}{\ca\!+J_X},$    $\quots{\k_o\bl t,x\br}{\ca\!+J_X},$ considered as $\k_o$-algebra.
      Suppose we get a system of equations $F(x,y)=0,$ whose solutions are $\k_o$-automorphisms of $R_X,R_Y,R_{X\!\times\! Y}.$ Add to
       this system the conditions $t_i\to t_i,$ i.e. $\Phi_X(t_i)=t_i$ and/or $\Phi_Y(t_i)=t_i$ and/or $C(t_i,o,o)=t_i.$ The enlarged system is still of implicit function type, and its solutions are now automorphisms of $\k$-algebras $R_X,R_Y,R_{X\!\times\! Y}.$ Then, having a formal (or approximate)   solution to this enlarged system of equations, we invoke the relevant Artin approximation of \S\ref{Sec.Background.AP}.
     \eee

\section{Strong Artin approximation and the P\l oski-version for groups $\cG$}\label{Sec.SAP.G.APP.G} Below $\cG$ is one of the group-actions $\cR,\cL,\cL\cR,\cC,\cK\circlearrowright\Maps.$ We define the $\cG$-approximation properties and then establish these for Nash/analytic/formal maps.

\subsection{Definitions}\label{Sec.GSAP.GAPP.definitions}
\bed ({\bf\SAPG})
The strong Artin approximation for $\cG$-action, \SAPG,   holds for a (Nash/analytic/formal) map $f\in \Maps$ if for any other map 
 $\tf\in \Maps$
 there exists a function $\N\stackrel{\be}{\to}\N$ ensuring: if $\tf\stackrel{mod\ \cm^{\be(d)}}{\equiv}g_d f,$ for some  $g_d\in \cG,$ then $\tf=gf$ for an element $g\in \cG$ that satisfies
 $g\stackrel{mod\ \cm^d}{\equiv} g_d. $
\eed
  Restate this via the filtration-topology. Take the orbit $\cG f$ and its filtration neighborhoods $\cG f+\cm^\bullet\cdot \RmX.$
   Suppose $\tf\in \cG f+\cm^d\cdot \RmX$ for each $d\in \N.$ Then $f\in \cap_d (\cG f+\cm^\bullet\cdot \RmX),$ the latter being the filtration closure
    of $\cG f\sseteq  \RmX.$ The  Strong Artin approximation ensures that the orbit is closed,  $\cap_d (\cG f+\cm^\bullet\cdot \RmX)=\cG f.$
 A remark: if $f$ is $\cG$-finitely determined, then the orbit is trivially closed, as $\cG f \spset \cG f+\cm^d\cdot \RmX$ for $d\gg1.$

\medskip 

Here is the explicit form of the  conditions, for the groups $\cL\cR$ and $\cK$:
\bei
\item \SAPLR.\quad  If $\Phi_{Y,d}(\tf(x)) \equiv  f( \Phi _{X,d}(x))\ mod\ \cm^{\be(d)}_X, $ for an element $(\Phi _{X,d},\Phi _{Y,d})\in Aut_X\times Aut_Y,$ then $\Phi_{Y }\circ\tf= f\circ \Phi _{X},$ for an element $(\Phi _{X },\Phi _{Y })\in Aut_X\times Aut_Y$
     that satisfies: $ \Phi _{X,d} \equiv  \Phi_X\ mod\ \cm^d_X $ and     $\Phi _{Y,d} \equiv \Phi_Y\ mod\ \cm^d_Y . $
\item \SAPK.
If $ \tf(\Phi_{X,d}(x)) \equiv  C_d(x,f(x)) \ mod\ \cm^{\be(d)}_X, $ for some $(\Phi _{X,d},C_d)\in \cK,$ then
 $ \tf(\Phi_{X }(x))= C (x,f(x)) ,$  for  an element $(\Phi _{X },C )\in \cK$
      that satisfies: $ \Phi _{X,d} \equiv  \Phi_X\ mod\ \cm^d_X  $ and    $C_d(x,y) \equiv C(x,y)\ mod(\cm_X,\cm_Y)^d. $
\eei
\bed({\bf\APPG})
The  parameterised approximation for $\cG$-action,  \APPG,    holds for a (Nash, resp. analytic)
 map $f\!\in\! \Maps$ if, given a formal solution $\tf\!=\!\hg f$ with $\hg\in \hat\cG,$ there exists a parameterized (Nash, resp. analytic) solution, $\tf\!=\!g_z f,$ that specializes to $\hg,$ i.e. $\hg\!=\!g_{\hz }.$
\eed
Explicitly (for $\cG=\cL\cR$): if $\hPhi_Y\circ \tf=f\circ\hPhi_X$ for some formal automorphisms $\hPhi_Y, \hPhi_X,$ then
 $\Phi_{Y,z}\circ \tf=f\circ\Phi_{X,z}$ for some (analytic/Nash) families of automorphisms with free parameters $z$, and moreover:
  $\hPhi_Y=\Phi_{Y,\hz(x,y)}$ and  $\hPhi_X=\Phi_{X,\hz(x,y)}$ for some   power series $\hz(x,y).$

\subsection{The approximations \SAPG and \APPG}
\bthe\label{Thm.SAP.G.and.APP.G} Let $R_X,R_Y$ be as in \S\ref{Sec.Background.Rings}.
\bee[\bf 1.]
\item Properties \SAPR, \SAPC, \SAPK hold for any (Nash/analytic/formal) map-germ $f\!\in\!\Maps.$
\\Properties \APPR, \APPC, \APPK hold for any (Nash/analytic) map-germ $f\in \Maps.$

\medskip

\item
Properties  \SAPL, \SAPLR hold  for any (Nash/analytic/formal)  map-germ $f\in \Maps.$
 In the formal case we assume that the base field   is a  $\aleph_0$-complete, see \S\ref{Sec.Background.SAP.nested}.
  In the analytic case we assume in addition: the map $f$ is of finite singularity type. 
%??? And liftability  
% when $\k$ is a field or $\k=\quots{\k_o\bl t\br}{\ca}$ with $\k_o$  any  field.
%\\
Properties \APPL, \APPLR hold for any  Nash   map-germ $f\in \Maps.$
 
 \eee
\ethe
\bpr We should resolve the condition $\tf=g f,$ $g\in \cG,$ while having an approximate (resp. formal) solution. We convert this condition into a system of equations of implicit function type. For the first statement we use the classical Strong Artin approximation and the P\l oski-Popescu approximation. For the second statement we use their nested versions. 

\bee[\bf 1.]
\item (We prove \SAPK and \APPK, the other versions are simpler.) We should resolve the condition $\tf(\Phi(x))=C(x,f(x))$ for unknowns $\Phi(x),C(x,y).$
 This is not an implicit function equation.      Present this condition in the form $\tf(\Phi (x)) -C(x,y)\in (y-f(x))\sset R_{X\!\times\! Y}.$  (This goes by Taylor expansion of $C(x,y)$ at $y=f(x)$.)
We get a system of equations of implicit function type:
\beq\label{Eq.inside.proof.IF.eq.to.resolve}
\tf(\Phi (x)) -C(x,y)= (y-f(x))\cdot a(x,y).
\eeq
Here $f,\tf$ are fixed, while $\Phi(x) , C(x,y),a(x,y)$ are unknowns. (See \S\ref{Sec.Background.AP.Remarks}.i.)

To this equation one adds also the condition $(Id,H)\in \cC,$ which means equation \eqref{Eq.group.C.def.equation}.
\bei
\item (\SAPK)
  Take the $\be$-function for the equations \eqref{Eq.group.C.def.equation} and \eqref{Eq.inside.proof.IF.eq.to.resolve}.
 Suppose $\tf \equiv g_d f\ mod\ \cm^{\be(d)}_X$ for some $g_d\in \cK,$ i.e.
 $\tf(\Phi_d(x)) \equiv C_d(x,f(x))\ mod\ \cm^{\be(d)}_X$ for some $(\Phi_d,C_d)\in \cK.$

 Then $(\Phi_d(x),C_d(x,y))$ is an approximate solution to \eqref{Eq.inside.proof.IF.eq.to.resolve} and \eqref{Eq.group.C.def.equation}.
  Invoke SAP, \S\ref{Sec.Background.SAP}, and remark \S\ref{Sec.Background.AP.Remarks}.iii. when $\k$ is a ring.
   We get the (Nash/analytic/formal) power series
 $\phi(x,y), C(x,y),a(x,y),$ satisfying   \eqref{Eq.inside.proof.IF.eq.to.resolve} and  \eqref{Eq.group.C.def.equation}, and with approximation properties:
\beq
\phi(x,y)\stackrel{mod(\cm_X,\cm_Y)^d}{\equiv}\Phi_d(x), \hspace{1cm} h(x,y)\stackrel{mod(\cm_X,\cm_Y)^d}{\equiv}C_d(x,y).
\eeq
  Finally define $\Phi_X(x):=\phi(x,f(x))$ to get  $\tf(\Phi(x))=C(x,f(x)).$

Here $(Id,C)\in \cC$ by construction, while the map $\Phi$ is invertible, being $\cm$-adic approximation to the automorphism $\Phi_d\in Aut_X.$

 \item (\APPK) Suppose $\tf\stackrel{\hat\cK}{\sim}f,$ i.e.  $\hf(\hPhi(x))=\hC(x,f(x))$ for some $(\hPhi,\hC)\in \widehat\cK.$
  Then equations \eqref{Eq.inside.proof.IF.eq.to.resolve},  \eqref{Eq.group.C.def.equation} have a formal solution.
  Invoke the APP, \S\ref{Sec.Background.APP}, to get a parameterized (Nash/analytic) solution, $\Phi_z(x,y), C_z(x,y),a_z(x,y)$ such that  $(\hPhi,\hC)$ is its specialization. Then  $(\Phi_z(x,f(x)), C_z(x,f(x)))$ is the needed family of elements in $\cK.$

\eei

\medskip

 \item (We prove \SAPLR and \APPLR, the other versions are  simpler.)
   We should resolve the condition $\Phi_Y(\tf(x))=f(\Phi_X(x))$ for the unknowns $(\Phi_X,\Phi_Y).$ This is not an implicit function equation. Present it in the form
   \beq\label{Eq.inside.proof.SAP.LR}
   \Phi_Y(y)-f(\Phi_X(x))\in (y-\tf(x))\sset R_{X\!\times\! Y}.
    \eeq
    (The two forms are equivalent by the Taylor expansion  $\Phi_Y(y)=\Phi_Y(\tf)+\cdots$.)

     We get the system of equations $\Phi_Y(y)-f(\Phi_X(x))=(y-\tf(x))\cdot a(x,y),$  for the unknowns $\Phi_X\in Aut_X,$ $\Phi_Y\in Aut_Y,$ $a\in R_{X\!\times\! Y}.$ (See Remark \S\ref{Sec.Background.AP.Remarks}.i.)
\bei
\item (\SAPLR) 
 Take the $\be$-function for the equations \eqref{Eq.group.C.def.equation} and \eqref{Eq.inside.proof.SAP.LR}.
 Suppose we have an approximate   solution,  
  \beq 
  \Phi_{Y,d}(y)-f(\Phi_{X,d}(x))\stackrel{mod(\cm_X,\cm_Y)^{\be(d)}}{\equiv}(y-\tf(x))\cdot a_d(x,y)\sset R_{X\!\times\! Y}.
  \eeq
     Here $\Phi_{X,d} {\in}\cR\ mod\ \cm^{\be(d)}_X$  and $\Phi_{Y,d}  {\in}\cL\ mod\ \cm^{\be(d)}_Y.$

Consider  $\Phi_{Y,d}(y),\Phi_{X,d}(x),a_d(x,y)$ as an approximate nested solution for the nest $\{y\}\sset \{x,y\}.$
     Invoke the nested SAP, \S\ref{Sec.Background.SAP.nested} (Nash in the Nash case, and formal in the analytic/formal cases), and remark \S\ref{Sec.Background.AP.Remarks}.iii. when $\k$ is a ring. We get the nested (Nash or formal) solution, $\phi_X(x,y), \Phi_Y(y),a(x,y)$ that approximates
    $\Phi_{Y,d}, \Phi_{X,d}.$ Finally  we specialize $\phi_X(x,y)\rightsquigarrow  \phi_X(x,f(x))=:\Phi_X(x).$
 Then $(\Phi_Y,\Phi_X)\in \cL\cR$ is the needed (Nash or formal) group-element.

Finally, in the analytic case (for maps of finite singularity type), apply the property \APLR, see Theorem 5.1 in \cite{Kerner.LRAP}.

\item (\APPLR) We have a formal solution of \eqref{Eq.inside.proof.SAP.LR}:
\beq
\hPhi_Y(y)-f(\hPhi_X(x))\in (y-\tf(x))\sset \quots{\k[\![x,y]\!]}{J_X+J_Y}, \quad\quad\quad
 \hPhi_Y\in \hat\cL,\quad \hPhi_X\in \hat\cR.
\eeq
  Consider it as a nested solution, for the nest $\{y\}\sset \{x,y\}.$
 Invoke the nested APP, \S\ref{Sec.Background.APP}  and remark \S\ref{Sec.Background.AP.Remarks}.iii when $\k$ is a ring. We get a  parameterized nested Nash solution $\Phi_{Y,z}(y),\phi_{X,z}(x,y).$
  Finally, specialize $\phi_{X,z}(x,y)\rightsquigarrow \phi_{X,z}(x,f(x))=:\Phi_{X,z}(x).$
 Then
 \beq
 \Phi_{Y,z}\circ\tf=f\circ \Phi_{X,z},\hspace{1cm}   \hPhi_X=\Phi_{X,\hz(x,y)},\hspace{1cm} and  \quad \hPhi_Y=\Phi_{Y,\hz(x,y)},
 \eeq
  for a specialization $\hz(x,y)\in \k[\![x,y]\!].$
  \epr
\eei
\eee

\subsection{A version of Theorem \ref{Thm.SAP.G.and.APP.G} for refined filtrations}\label{Sec.SAP.G.APP.G.remarks}
  The proved versions of Artin approximation, \SAPG and \APPG, are
  for the filtrations $\cm^\bullet_X,$ $\cm^\bullet_Y.$
  As is explained in Remark \ref{Rem.Introduction}, we cannot strengthen these versions to the filtrations $I^\bullet_X,\cm^\bullet_Y,$
   when $\sqrt{I_X}\ssetneq \cm_X.$
 Yet, we give a refinement to the filtrations  $\ca_X\cdot \cm^\bullet_X,$ $\ca_Y\cdot \cm^\bullet_Y .$
   Here  the ideals $\ca_X\sset R_X,$ $\ca_Y\sset R_Y$ are not   $\cm_X,\cm_Y$-primary, i.e. $\sqrt{\ca_X}\ssetneq \cm_X$ and  $\sqrt{\ca_Y}\ssetneq \cm_Y.$
  The argument is standard, see \cite{Belitskii.Boix.Kerner}. We give the details only for the group $\cL\cR,$ the other cases are simpler.
   Fix a Nash map $X\stackrel{f}{\to}Y$ and ideals  $\ca_X\sset R_X,$ $\ca_Y\sset R_Y.$
 \bee[\bf i.]
  \item ({\pmb\SAPLR$.\ca_X\ca_Y$})  There
exists a function $\N\stackrel{\be}{\to}\N$ ensuring: if
$\Phi_{Y,d}(\tf(x))\stackrel{mod\ \cm^{\be(d)}_X}{\equiv} f( \Phi _{X,d}(x)) $
for some $\tf\in \Maps$ and  elements satisfying $\Phi _{X,d}(x)-x\in \ca_X,$ $\Phi _{Y,d}(y)-y\in \ca_Y,$
  then $\Phi_Y\circ \tf= f\circ \Phi_X,$ where the elements  $\Phi_Y,\Phi_X$ satisfy
$\Phi _X \stackrel{mod\ \ca_X\cdot\cm^{d}_X}{\equiv}\Phi _{X,d}   ,$
 $\Phi _Y\stackrel{mod\ \ca_Y\cdot\cm^{d}_Y}{\equiv}\Phi _{Y,d}.$

\bpr Fix some generators of ideals, $\ca_X=R_X\bl \{q^X_\bullet\}\br$ and $\ca_Y=R_Y\bl \{q^Y_\bullet\}\br.$ We are looking for automorphisms of the form $\Phi_X(x)=x+\sum q^X_j c^X_j,$ $\Phi_Y(y)=y+\sum q^Y_i c^Y_i,$ with the unknowns $c^X_\bullet,c^Y_\bullet.$
  As in the proof of theorem \ref{Thm.SAP.G.and.APP.G} we convert the condition $\Phi_Y\circ \tf=f\circ\Phi_X$ into a system of implicit function equations on $c^X_\bullet,c^Y_\bullet,$ see \eqref{Eq.inside.proof.SAP.LR}.
  These equations have  an approximate (nested) solution $c^X_{\bullet,d}(x),c^Y_{\bullet,d}(y).$ Apply the nested SAP, and then specialize $y\rightsquigarrow f(x),$ to get
   the needed coefficients $c^X_\bullet(x),c^Y_\bullet(y).$
\epr

\item ({\pmb \APPLR.$\ca_X\ca_Y$}) Suppose $\hPhi_Y\circ \tf=f\circ\hPhi_X$ for some formal automorphisms $\hPhi_Y, \hPhi_X,$ satisfying
 $\hPhi_Y(y)-y\in \ca_Y\cdot \hR_Y,$ $\hPhi_X(x)-x\in \ca_X\cdot \hR_X.$ Then  there exists a parameterized solution
   $\Phi_{Y,z}\circ \tf=f\circ\Phi_{X,z} $ satisfying :
   \[
   \Phi_{Y,z}(y)-y\in \ca_Y,\quad     \Phi_{X,z}(x)-x\in \ca_X,\quad
   \hPhi_Y=\Phi_{Y,\hz(x,y)}, \quad \hPhi_X=\Phi_{X,\hz(x,y)}
   \]
     for some   power series $\hz(x,y)\in \k[\![x,y]\!].$

    The proof is the same as in part i.
\item In the same way one gets the versions of \SAPR, \SAPL, \SAPC, \SAPK for the filtrations  $\ca_X\cdot \cm^\bullet_X,$ $\ca_Y\cdot \cm^\bullet_Y,$ and the versions of \APPR, \APPL, \APPC, \APPK for $\ca_X,\ca_Y.$   We omit the details.
     \eee

\subsection{A version of Theorem \ref{Thm.SAP.G.and.APP.G} for unipotent subgroups $\pmb{\cG^{(j)}<\cG}$}\label{Sec.SAP.G.APP.G.unipotent}
 Fix an ideal $I\sset R_X,$ and take the filtration $\RmX\spset I\cdot \RmX\spset I^2\cdot \RmX\spset\cdots.$ 
  This defines the filtration on the space of maps, $M^d:=\Maps\cap I^d\cdot \RmX.$ For $Y=(\k^m,o)$ we have $\Maps=\cm\cdot \RmX,$ 
   therefore $M^d=I^d\cdot \RmX.$ For non-$\k$-smooth targets the subsets $M^d\sset I^d\cdot \RmX$ are neither additive, nor $\k$-multiplicative. 
  
 This filtration defines the subgroups of the groups $\cL,\cR,\cC$ for $j\ge0$ by
 \beq\label{Def.filtered.groups}
 \cG^{(j)}:=\{g\in \cG|\ g\cdot M^d=M^d  \text{ and }Id=[g]\circlearrowright \quot{M^d+I^{d+j}\cdot \RmX }{I^{d+j}\cdot \RmX } \quad \forall\ d\ge1 \}.
 \eeq
Thus we get the filtration by subgroups $\cG\ge \cG^{(0)}\ge \cG^{(1)}\ge\cdots.$ See \cite{Kerner.Group.Orbits} and \cite{BGK.20} for   detail.

\bex\label{Ex.Filtered.Groups}
 Suppose  the target is $\k$-smooth, $Y=(\k^m,o).$ Thus $\Maps=\cm_X\cdot\RmX.$
  Fix some generators $I=R_X\cdot \{q_i\}.$ 
 \bee[\!\!\!$\bullet$]
\item $\cR^{(j)} =\{\Phi_X\in Aut_X|\ q_i(\Phi_X(x))-q_i\in I^{j+1}\quad \forall i \}.$ 
\\For $I=\cm_X$ this simplifies to  $\cR^{(j)} =\{\Phi_X\in Aut_X|\  \Phi_X(x)-x\in \cm_X^{j+1} \},$  
i.e. automorphisms of $X$ that are $Id$ modulo $\cm^{j+1}_X.$

\item $\cL^{(j)} =\{\Phi_Y\in Aut_Y|\ \Phi_Y(y)=y+Q(y),\ \text{where } Q(I^d\cdot \RmX)\sseteq I^{d+j}\cdot \RmX\quad \forall\ d\ge1 \}.$
\\Thus $\cL^{(j)}\ge\{\Phi_Y\in Aut_Y|\ \Phi_Y(y)-y\in (y)^{j+1}\}.$ For $I=(x)$ and    $\k\spseteq\Q$  the two groups coincide.
   See example 3.2 in \cite{Kerner.Group.Orbits}.
\item 
Similarly  $\cC^{(j)}\! =\{(Id_X,C)\!\in\! Aut_{X\!\times\! Y}|\ C(x,y)-y\!=\!Q(x,y), \ \text{where } Q(x,I^d\!\cdot\! \RmX)\!\in  \! I^{d+j}\!\cdot\! \RmX \}.$
\\Thus $\cC^{(j)}\ge \{(Id_X,C)\in Aut_{X\!\times\! Y}|\ C(x,y)-y\in (x+y)^{j+1}\}.$
\eee
\eex
When the target is singular, $Y\ssetneq(\k^m,o),$ the groups $\cR^{(j)}, \cL^{(j)}, \cC^{(j)}$ are more complicated. 

For the composite groups we define $(\cL\cR)^{(j)}:=\cR^{(j)}\times \cL^{(j)}$ and  $\cK^{(j)}:=\cC^{(j)}\rtimes \cR^{(j)}.$
 
 \bcor
 In the assumptions of Theorem \ref{Thm.SAP.G.and.APP.G} the properties SAP.$\cG^{(j)}$ and APP.$\cG^{(j)}$ hold for $\Maps.$
 \ecor
\bpr
In view of the proof of theorem \ref{Thm.SAP.G.and.APP.G} it is enough to verify: the additional condition $g\in \cG^{(j)}$ is an implicit function equation. The verification is separate for each group.
\bee[\hspace{-0.3cm}$\bullet$]
\item 
$\cC^{(j)} =\{(Id_X,C(x,y))\in Aut_{X\!\times\! Y}|\ C(x,f(x))-f(x)\in I^{d+j}\cdot\RmX\quad \forall\ d\in \N,\ \forall\ f\in M^d\}.$ This condition on $C$ can be rewritten
 as $C(x,y)-y\in (y-f(x))+I^{d+j}\cdot\RmXY,$  by Taylor expansion. Therefore we can present this group as
 \[
 \cC^{(j)} =\{(Id_X,C(x,y))\in Aut_{X\!\times\! Y}|\ C(x,y)-y\in \capl_{d\in \N}\capl_{f\in M^d} \left((y-f(x))+I^{d+j}\cdot\RmXY\right)\}.
 \]
The subset $P:=\cap_{d\in \N}\cap_{f\in M^d} \left((y-f(x))+I^{d+j}\cdot\RmXY\right)\sset \RmXY$ is an $R_{X\!\times\! Y}$-submodule, being the intersection of submodules. Therefore $P$ is finitely generated. And thus the condition $C(x,y)-y\in  P$ is presentable as a (finite) system of equations of implicit function type.

\item
$\cL^{(j)} =\{\Phi_Y\in Aut_Y|\ \Phi_Y(f(x))-f(x)\in I^{d+j}\cdot\RmX\quad  \forall\ d\in \N,\ \forall\ f\in M^d\}.$ 
 As in the case $\cC^{(j)} $ this group  can be presented  in the form 
$
 \cL^{(j)} =\{\Phi_Y\in Aut_Y|\ \Phi_Y(y)-y\in P\}.
$ 
 Here $P\sset \RmXY$ is the same (finitely generated) $R_{X\!\times\! Y}$-submodule.
 Thus the condition $\Phi_Y(y)-y\in  P$ is presentable as a (finite) system of equations of implicit function type. But now we look for a nested solution, i.e. $\Phi_Y$ depends on $y$ only. 

\item 
$\cR^{(j)} =\{\Phi_X\in Aut_X|\  f(\Phi_X(x))-f(x)\in I^{d+j}\cdot\RmX\quad  \forall\ d\in \N,\ \forall\ f\in M^d\}.$ 
 The left hand side is a vector of power series, with a given $f$ and unknown $\Phi_X.$ These left hand sides generate the submodule
  $P:=\sum_{d\in \N}\sum_{f\in M^d}( f(\Phi_X(x))-f(x))\sseteq \RmX[\![\Phi_X]\!],$ resp. $\RmX\{\Phi_X\}, \RmX\bl\Phi_X\br.$
   This submodule is finitely generated. Therefore the subgroup $\cR^{(j)}<\cR$ is defined by (a finite number of) equations of implicity function type. 
\eee
For the composite subgroups, $\cL\cR^{(j)}<\cL\cR,$ $\cK^{(j)}<\cK,$ we collect the conditions of their factors, getting again implicit function equations.  
\epr

\section{Artin approximation for quivers of maps}\label{Sec.Quivers}
Below $\Ga$ is   a finite directed graph with at least two vertices.

 We define the approximation problem of \S\ref{Sec.Intro.AP.for.quivers} and establish this approximation for directed rooted trees.

\subsection{Quivers of maps}\label{Sec.Quivers.Definition}
As in \S\ref{Sec.Intro.AP.for.quivers}, the vertices of $\Ga$ are (formal, resp. analytic, Nash) scheme-germs, $\{X_v=Spec(R_v)\}_v.$  The arrows are morphisms,
 $\{X_v \stackrel{f_{wv}}{\to} X_w\}_{wv}.$
    Call this data ``a quiver of maps on the graph $\Ga$".
\\\parbox{12cm}{\bed\label{Def.Morphism.of.Quivers}
    A morphism of two such quivers on the same graph (and over the same base $\k$), $(\Ga,\{\tX_\bullet\},\{\tf_{wv}\}) \to (\Ga,\{X_\bullet\},\{f_{wv}\}),$ is a collection of maps $ \tX_v\stackrel{\Phi_v}{\to} X_v,$ making all the vertical rectangles of the diagram to commute.
\eed}\hspace{1cm} %\beq\label{Eq.Quiver.AP.Diagram}
$\bM
 \underset{\nearrow}{\stackrel{\searrow}{\to}}\quad
  \tX_v\stackrel{\tf_{wv}}{\to} \tX_w\quad
 \underset{\searrow}{\stackrel{\nearrow}{\to}}
\\ \Phi_v\downarrow\quad\quad  \downarrow \Phi_w
\\
 \underset{\nearrow}{\stackrel{\searrow}{\to}}\quad
  X_v\stackrel{f_{wv}}{\to} X_w\quad
 \underset{\searrow}{\stackrel{\nearrow}{\to}}
\eM$ %\eeq

 This commutativity of the vertical rectangles gives the system of equations:
  $\{\Phi_w\!\circ\! \tf_{wv}=f_{wv}\!\circ\! \Phi_v\}_{v,w},$
   with unknowns $\{\Phi_\bullet\}.$
 Each  $\Phi_v$ corresponds to a (local) homomorphism of   local $\k$-algebras, $ R_v\stackrel{\Phi^\sharp_v}{\to} \tR_v,$
  $x_v\!\to\! \Phi_v(x_v).$
  Fix some generators of the maximal ideals of all the local rings $R_\bullet $, then the morphisms $ \Phi^\sharp_\bullet $ are presented as  vectors of (analytic/algebraic/formal) power series.

\beR\label{Rem.Additional.Equations.On.Quivers}
One can impose various additional conditions on the maps $ \Phi_\bullet ,$ e.g.:
\bee[\bf i.]
\item Some of the maps $\Phi_v$ vanish at the origin in a prescribed manner, i.e. $\Phi^\sharp_v(\cm_{X_v})\sseteq \tca_v\sset \tR_v.$
 (The ideals $\tca_\bullet$ are prescribed.)
 \item Some of the maps $\Phi_v$ send subgerms $\tZ_v\sset \tX_v$ to subgerms $Z_v\sset X_v,$ i.e. $\Phi^\sharp_v(I_{Z_v})\sseteq I_{\tZ_v}\sset \tR_v.$
\item Some of $\Phi_v$ are isomorphisms. 

When all $\Phi_v$ are isomorphisms, we get the ``$\cA$-equivalent of quivers".
\item (Assuming $X_v=\tX_v$ for all vertices.) Some $\Phi_v$'s are identities up to ``higher-order-terms", i.e. $\Phi_v(x_v)-x_v\in \ca_v$ for some ideals $\ca_v\sset R_v.$

    In the particular case $\ca_v=0$ this means $\Phi_v=Id_{X_v}.$
  Then $\Phi_v$ is just removed from the list of unknowns.
\eee
Observe that all these additional conditions on $ \Phi_\bullet $ are presentable as equations of implicit function type, see \S\ref{Sec.Background.AP.Remarks}.ii.
\eeR
\bex\label{Ex.Quivers.R,L,LR,cases}
Take the simplest quiver $X\stackrel{f}{\to}Y,$ the corresponding equation is $\Phi_Y\circ \tf=f\circ\Phi_X.$
\bei
\item
Imposing the condition
 ``$\Phi_X,\Phi_Y$ are isomorphisms"    gives the left-right equivalence of maps $f,\tf,$ see \S\ref{Sec.Background.R.LR.equiv}.
\item Imposing in addition $\Phi_X=Id_X$ (resp. $\Phi_Y=Id_Y$)  gives the left (resp. the right) equivalence of maps $f,\tf.$
\item Take $f=0$ and impose no invertibility condition on $\Phi_Y.$ One gets the approximation problem for the equation
 $\Phi_Y\circ\tf=0,$ i.e. the inverse Artin problem for $\tf$.
\item Fix some subgerms $Z_X\sset X$ and $Z_Y\sset Y,$ and impose the conditions $\Phi_X(Z_X)=Z_X,$  $\Phi_Y(Z_Y)=Z_Y.$
 We get the $\cL\cR$-equivalence relative to these subgerms.
\eei
\eex
\bed\label{Def.Gamma.AP} (Approximation properties for quivers)
\bee
\item {\pmb \APGam}. The Artin approximation holds for the quiver $(\Ga,\{X_v\},\{f_{wv}\})$ if any formal morphism,
 $ (\Ga,\{\tX_v\},\{\tf_{wv}\})\stackrel{\{\hPhi_v\}_v}{\to }(\Ga,\{X_v\},\{f_{wv}\}),$ is approximated by Nash/analytic morphisms.
 \item
The strong Artin approximation for quivers, {\bf\SAPGam}, and the P\l oski-type approximation, {\bf\APGamP}, are defined similarly.
 \eee
\eed
  Explicitly (for  \APGam): given any quiver $(\Ga,\{\tX_v\},\{\tf_{wv}\})$
 (over the same base ring $\k$),
  any set of formal morphisms $\{\hat{\tX}_i\stackrel{\hPhi_v} {\to} \hX_v\}$
 satisfying the equations $\{\Phi_w\circ \tf_{wv}=f_{wv}\circ \Phi_v\}_{wv},$ is approximated by sets of ordinary solutions
  $\{ \tX_v\stackrel{\Phi_v}{\to }X_v\}.$

Depending on the context one adds the conditions of remark \ref{Rem.Additional.Equations.On.Quivers}.

\beR\label{Rem.Gamma.AP}
 From the very beginning we exclude graphs with loops.
 Any loop will produce the Schr\"oder-type equation   $f\circ \Phi_X=\Phi_X\circ \tf.$ Imposing the condition ``$\Phi_X$ is an automorphism", we would get the conjugation equivalence, $\tf=\Phi_X\circ f \circ\Phi^{-1}_X.$ Recall that in this case the formal and $\C$-analytic classifications differ essentially even in the one-dimensional case.  See Siegel's linearization theorem, \cite[pg.114]{Milnor},  \cite[pg.356]{AGV}.

  For the same reasons we restrict to the case: $\Ga$ is a tree.

% The notion of morphism of quivers, definition \ref{Def.Morphism.of.Quivers}, extends the $\cL\cR$-equivalence of maps. The $\cK$-equivalence of %maps is extended  similarly to the $\cK$-equivalence of quivers. Then  definition \ref{Def.Gamma.AP} gives the  \APGamK  property.
 %    And proposition \ref{Thm.Gamma-AP.for.k<x>} holds also for  \APGamK.
\eeR

\subsection{Properties \APPGam and \SAPGam for   directed rooted trees}\label{sec.Gamma.AP.for.rooted.trees}
Below $\Ga$ is a directed rooted tree. Namely, every vertex $v\in \Ga$ is connected by a unique directed  path to a particular vertex, $v\rightsquigarrow o\in \Ga,$ the later is called the root of $\Ga.$

\bthe\label{Thm.Gamma-AP.for.k<x>}
\bee
\item  The  property \APPGam holds in the $\k$-Nash case, i.e. $\{R_v=\quots{\k\bl x_v\br}{J_v}\}_v.$

\item   The property \SAPGam holds in the Nash case  and (when the base field is   $\aleph_0$-complete) in the formal case.  (See \S\ref{Sec.Background.SAP.nested})
   \eee
\ethe
\bpr
The proof is an extension of the proof of Theorem \ref{Thm.SAP.G.and.APP.G}. We convert the conditions to resolve into a system of equations of implicit function type. Their solutions should depend only on particular types of variables. Transform this into a nested approximation problem. Invoke the nested versions of SAP and APP. Then refine the obtained (nested) solutions iteratively, to get solutions in separated variables.
\bee[\hspace{-0.3cm}\bf Step 1.]
\item

Given two such trees  we get the diagram of definition \ref{Def.Morphism.of.Quivers}, with  the system of equations:
\beq\label{Eq.inside.proof.system.to.resolve}
 \Phi_w\circ \tf_{wv}=f_{wv}\circ \Phi_v  \hspace{3cm} \text{ for each edge }[v\to w]\sset \Ga,
%, \quad\quad\quad  \{x_w=f_{wv}(x_w)\}, \quad\quad\quad  \{\tx_w=\tf_{wv}(\tx_v)\},
\eeq
and possibly additional equations of Remark \ref{Rem.Additional.Equations.On.Quivers}.
 Explicitly, $\Phi_w( \tf_{wv}(\tx_v))=f_{wv}(\Phi_v(\tx_v)) $ where:
 \bei
 \item  the multivariable   $x_v$ (for each $v$)  denotes the generators of the   ideal $(x_v)\sset R_v;$
\item the maps $\{f_{wv},\tf_{wv}\}$ are prescribed, while  $\{\Phi_v\}_v$ are unknown homomorphisms.
\eei
The condition \eqref{Eq.inside.proof.system.to.resolve} is not an implicit function equation. As in the proof of Theorem \ref{Thm.SAP.G.and.APP.G} we rewrite this condition in the form
 $\Phi_w(\tx_w)-f_{wv}(\Phi_v(\tx_v))=(\tx_w-\tf_{wv}(\tx_v))\cdot a_{wv} ,$ with additional unknowns $\{a_{wv}\}.$

 To this system  we add the conditions $\{\Phi_v\circlearrowright R_v\},$ i.e. $\{\Phi_v(J_v)\sseteq J_v\},$
  and also the conditions of  Remark \ref{Rem.Additional.Equations.On.Quivers}.
 Altogether we get a system of implicit function equations.

By our assumption this system admits a formal solution for part 1 and
an approximate solution for part 2. Explicitly (for part 2): $\{\Phi_{v,d}\}_v$
  and $\{a_{wv,d}\}_{wv},$ satisfying:
  \beq
  \Phi_{w,d}(\tx_w))-f_{wv}(\Phi_{v,d}(\tx_v))\stackrel{mod\ (\tx_v^{\be(d)})}{\equiv} (\tx_w-\tf_{wv}(x_v))\cdot a_{wv,d}
   \hfill \text{ for each edge }[v\to w]\sset \Ga.
  \eeq
We should approximate the elements $\{\Phi_{v,d}\}_v$ by an ordinary (Nash, resp. formal) solution. (For part 1 we should present $\{\hPhi_v\}$ as specializations of Nash automorphisms.)

We will use the nested P\l oski approximation \S\ref{Sec.Background.APP}, resp. the nested
Strong Artin approximation, \S\ref{Sec.Background.SAP.nested}. Thus we should define  the nested structure on the set of variables.

Define the grade of a vertex $v\in \Ga$ as the length of the (unique) directed path to the root, $v\rightsquigarrow o.$ Thus the grade of the root is 0. Denote by $\Ga_{\le j}\sset \Ga$ the subgraph of all the vertices (and their edges) of grades$\le j.$ Thus $\Ga_{\le 0}=\{o\},$ while
 $\Ga_{\le 1}$ is a ``star" graph, $K_{1,n}$.

Define the nested structure on the set of variables $\{\tx_\bullet\}$ by the chain of subgraphs $\Ga_{\le 0}\sset \Ga_{\le1}\sset\cdots\sset\Ga.$
 Observe: the subgraph $\Ga_{\le j}\smin \Ga_{\le j-1}$ is totally disconnected. (I.e. if $v,w\in \Ga_{\le j}\smin \Ga_{\le j-1}$ then $v,w$ are not connected by an edge.)

The approximate solution $\{\Phi_{v,d}\}$ (for part 2), resp. the formal solution $\{\hPhi_{v}\}$ (for part 1) are nested for this nest.

\item (Below we treat part 2. The proof of part 1 is almost the same.)
We have the approximate nested solution,   $\{\Phi_{v,d}(\tx_v)\}_v.$ Invoke   the nested SAP, \S\ref{Sec.Background.SAP.nested}, to get the nested  approximation, i.e. Nash/formal power series  $\{\Psi_v(\{\tx_{(\bullet)}\}_{gr(\bullet)\le gr(v)})\}_v.$
  (For the case ``$\k$ is not a field" see Remark \ref{Sec.Background.AP.Remarks}.iii.)
These power series are ``non-pure", i.e.
     $ \Psi_v(\dots)$ depends not only on $\tx_v.$
        Thus they do not define  morphisms of the scheme-germs $\{X_v\}_v.$

     We convert these into pure power series, using specializations.

 By our construction, for each edge of the graph     one has the inclusion:
\beq\label{Eq.inside.proof.Inclusion}
\quad
[v\to w]\sset \Ga \hspace{2cm}
\Psi_w(\{\tx_\bullet\}_{gr(\bullet)\le gr(w)})-f_{wv}(\Psi_v(\{\tx_\bullet\}_{gr(\bullet)\le gr(v)}))\in
   (\tx_w-\tf_{wv}(\tx_v)).
\eeq
Both parts are power series. The variables $\{\tx_\bullet\}_{\bullet\neq v,w}$ are free parameters. Therefore this particular condition
 $[v\to w]$ remains true under {\em any} specialization of these parameters.

We will specialize the parameters iteratively, so that the functions $\Psi_v (\{\tx_\bullet\})$ are transformed into ``pure" solutions,
    $\{x_v(\tx_v)\}_v.$ (And the equations satisfied by $\Psi_v$ still hold.)
Observe that the solution for the root of $\Ga$ is already pure,   $\Psi_{(o)}(\{\tx_\bullet\}_{gr(\bullet)\le gr( o)})= \Psi_{(o)}(\tx_o).$

\

\underline{The iterative step.} Suppose the solutions $\Psi_v ( \tx_\bullet )$ are pure for vertices of grade $gr(v)\le j,$ while  for $gr(v)> j$ all the functions
 $\Psi_v ( \tx_\bullet )$ depend on $\tx_\bullet$ with $gr \bullet \ge j$ only. (Initially this holds for the root, $j=0.$)
 We specialize certain parameters to pass from $j$ to $j+1.$

 Fix any vertex $w$ of $grade=j$ and any $v$ of $grade>j$ that is connected to $w$ by a (unique) directed path, $v\rightsquigarrow w.$ Compose the maps along this path, $f_{w\cdots v}:=f_{wv_1}\circ f_{v_1v_2}\circ\cdots\circ f_{v_n v}.$
  Specialize the power series $\Psi_v ( \tx_\bullet )$ to $\Psi_v ( \tx_\bullet )|_{\tx_w=\tf_{w\cdots v}(\tx_v)}.$
   We claim: this specialization preserves all the conditions \eqref{Eq.inside.proof.Inclusion}.
   For the condition of the edge $[v\to w]$ this is verified by Taylor expansion of $\Psi_v(\tx).$
    For all the other edges the preservation is trivial.

Make such substitutions (repeatedly) for all vertices of the complement $\Ga_{\le j+1}\smin \Ga_{\le j}.$ At each step we replace all the power series
 $\Psi_v ( \tx_\bullet )$ by their specializations. This does not change the power series $\Psi_v ( \tx_\bullet )$ for $gr(v)\le j.$

By now, for each $v\in \Ga_{\le j+1}$ the power series $\Psi_v ( \tx_\bullet )$ depends only on $\tx_\bullet$ with $gr(\bullet)=j+1.$
 It is not yet pure. Finally, we specialize $\Psi_v ( \tx_\bullet )\to  \Psi_v ( \tx_\bullet )|_{\tx_{(\bullet\neq v)}=0},$
  to get a pure power series, $\Psi_v ( \tx_v ).$ Do this for all the descendants of $v$, i.e.:
    $\Psi_{(v')}( \tx_\bullet )\to  \Psi_{(v')}( \tx_\bullet )|_{\tx_{(\bullet\neq v)}=0},$ whenever $[v'\rightsquigarrow v].$
 Again, this preserves all the conditions of \eqref{Eq.inside.proof.Inclusion}.

\

The specialized power series satisfy now:
 $\Psi_v ( \tx_\bullet )$ are pure for $gr(v)\le j+1,$ while  for $gr(v)> j+1$ all the functions
 $\Psi_v ( \tx_\bullet )$ depend on $\tx_\bullet$ with $gr (\bullet) \ge j+1$ only.

  Repeat this iterative step. By the end of these iterations the non-pure solution $\Psi_\bullet$ is converted into a (pure) collection of homomorphisms, $\Phi_V:R_v\to \tR_v.$
\epr
\eee
\bex
Take a (finite) multi-germ $\amalg X_i,$ and a map $\amalg X_i\stackrel{\amalg f_i}{\to}Y.$ The corresponding quiver in \eqref{Eq.quiver.multigerms} is a directed rooted tree.
\bee[\hspace{-0.2cm}$\bullet$]
\item
 Suppose two such Nash maps,  $\amalg X_i\stackrel{\amalg f_i,\ \amalg \tf_i}{\to}Y,$ are formally $\cL\cR$-equivalent.
  Namely, $\{\hPhi_Y\circ \tf_i=f_i\circ\hPhi_{X_i}\}_i,$ for some formal automorphisms $\hPhi_{X_i}\circlearrowright \hX_i,$
   $\hPhi_{Y }\circlearrowright \hY .$ Then there exists a parameterized Nash equivalence,
    $\{\Phi_{Y ,z}\circ \tf_i=f_i\circ\Phi_{X_i,z}\}_i,$ satisfying:
          $\hPhi_{X_i}=\Phi_{X_i,\hz}$ and $\hPhi_{Y }=\Phi_{Y ,\hz}$ for some (formal) specialization of the parameter,
    $\hz.$

    In particular, the formal equivalence, via  $\{\hPhi_Y,\hPhi_{X_i} \}_i,$ is approximated by Nash equivalence,  $\{\Phi_{Y },\Phi_{X_i} \}_i,$ in the $\cm_X,\cm_Y$-adic topology.
\item
For a given (Nash or formal) map $\amalg X_i\stackrel{\amalg f_i}{\to}Y$ there exists a function $\N\stackrel{\be}{\to}\N$ ensuring:
 if the maps $\amalg f_i,$ $\amalg \tf_i$ are approximately $\cL\cR$-equivalent,  i.e. $\{\Phi_{Y ,d}\circ \tf_i\stackrel{\cm_X^{\be(d)}}{\equiv}f_i\circ\Phi_{X_i,d}\}_i,$  then they are $\cL\cR$-equivalent, i.e.
   $\{\Phi_{Y ,z}\circ \tf_i=f_i\circ\Phi_{X_i,z}\}_i,$ where moreover  $\Phi_{X_i}\stackrel{\cm_X^d}{\equiv}\Phi_{X_i,d}$ and $\Phi_{Y }\stackrel{\cm_Y^d}{\equiv}\Phi_{Y ,d}.$
\eee

\eex

\section{Artin approximation with base change}\label{Sec.AP.Gamma.with.parameters}

Definitions \ref{Def.Morphism.of.Quivers} and \ref{Def.Gamma.AP} address   (quivers of) maps over the fixed base  $Spec(\k).$ The maps
  $R_v\stackrel{\Phi_v}{\to} \tR_v$ were homomorphisms of $\k$-algebras. For various applications, e.g. when studying (induced) unfoldings/deformations, it  is important to allow also morphisms of the base ring $\k.$
 We define the corresponding approximation property, and deduce it for directed rooted trees.

\subsection{Unfoldings of maps} Let $\k$ be a local ring.
 The cases relevant for us are: $(\k,R_X,R_Y)$ as in \S\ref{Sec.Background.Schemes.Maps}.ii.
 Consider the map $X\stackrel{f}{\to}Y$  as an unfolding of its central fibre,
 \beq
Spec(R_X\otimes\quots{\k}{\cm_\k})=: X_o\stackrel{f_o:=f\otimes \quots{\k}{\cm_\k}}{\longrightarrow}Y_o:= Spec(R_Y\otimes\quots{\k}{\cm_\k}).
 \eeq
 \parbox{12.2cm}{ This is a family of maps over the base $Spec(\k).$ Another unfolding, $\tX\stackrel{\tf}{\to}\tY,$
is called $\cL\cR$-induced from $f$ if there is a commutative diagram of maps.  (See e.g. \cite{Mond-Montaldi}.)
 In some cases one insists that the central fibres coincide, i.e. $\Phi_{X_o}=Id_{X_o}$ and $\Phi_{Y_o}=Id_{Y_o}.$
}
 \quad\quad$\bM \hspace{0.5cm}\tX\quad\stackrel{\Phi_X}{\longrightarrow} \quad X\hspace{0.5cm}\\
  \pi_\tX \searrow \tf \hspace{0.7cm} \pi_X  \searrow f\\
 \hspace{1cm} \tY \stackrel{\Phi_Y}{\longrightarrow}\hspace{0.8cm}  Y\\
\hspace{0.7cm} \swarrow\pi_\tY\hspace{0.9cm} \swarrow \pi_Y\\
 Spec(\tilde\k)\stackrel{\Phi_\k}{\longrightarrow}Spec(\k)\eM \put(-95,5){\Big\downarrow}\put(-35,5){\Big\downarrow}$

Accordingly we write $f$ as $f_t(x),$  resp. $\tf$ as $\tf_\tt(\tx),$ and $\Phi_X$ as $\Phi_{X,\tt} ,$ and so on.
 Then the diagram gives the equation $\Phi_{Y,\tt}(\tf_\tt(\tx))=f_{\Phi_\k(\tt)}(\Phi_{X,\tt}(\tx)).$
  Here $\Phi_\k(\tt),\Phi_{X,\tt},\Phi_{Y,\tt}$ are unknowns.

More generally, for an equivalence group  $\cG,$ one wants to verify whether an unfolding $\tf_\tt$ is $\cG$-induced from $f_t.$
 The equation to resolve is $\tf_\tt=g_\tt(f_{\Phi_\k(\tt)})$ with the unknowns $g_\tt,\Phi_\k.$
 E.g. for $\cG=\cK$ the equation is" $\tf_\tt(\tx)=C_{\Phi_\k(\tt)}(x,f_{\Phi_\k(\tt)} ( \Phi_{X,\tt}(\tx)).$

 \subsection{Unfoldings of quivers}

More generally, for $\k$ a local ring, one considers quivers of families of (formal/analytic/Nash) map-germs over $Spec(\k),$ i.e. $\{X_{t, v}\stackrel{f_{t,wv}}{\to}X_{ t,w}\}_{vw},$ with $R_{t,  v }=\quots{\k\{x_v\}}{J_{ v}},$ resp. $\quots{\k\bl x_v\br}{J_{ v}}.$
 We want to $\cG$-induce a quiver deformation, i.e. to resolve the system of equations $\{\tf_{\tt,v}=g_{\tt,v}(f_{\Phi_\k(\tt)})\}_v$ with the unknowns $\{g_{\tt,v}\}_v,\Phi_\k.$

 \bed{\bf AP.$\pmb{ (\Ga,t)} $}
  The   Artin approximation with parameters holds for the group $\cG$ acting on the quiver-family $(\Ga,\{X_{t,v}\},\{f_{t,wv}\})$ if for any quiver-family  $(\Ga,\{\tX_{\tt,v}\}_v,\{\tf_{\tt,wv }\}_{wv},$
  any formal solution $\{\tf_{\tt,v}=\hg_{\tt,v}(f_{\hPhi_\k(\tt)})\}_v$  is approximated by ordinary (analytic/Nash) solutions.
\eed

The properties SAP.${(\Ga,t)}$, APP.${(\Ga,t)}$ are defined similarly.

Below $\cG=\cL\cR$ and we write the definition explicitly:  any
  quiver of formal morphisms $\{ \hat\tX_{t,v}\stackrel{\hPhi_v}{\to} \hX_{t,v}\},$ $Spec(\hat\tk)\stackrel{\hPhi_\k}{\to} Spec(\hk),$
 solving the equations $\{\Phi_w\circ \tf_{t,wv}=f_{\Phi_\k(\k),wv}\circ \Phi_v\}_{v,w},$ is approximated by a quiver of ordinary   morphisms
  $\{ X_{t,v}\stackrel{\Phi_v}{\to} \tX_{\tt,v}\}_v,$ $Spec(\k)\stackrel{\Phi_\k}{\to} Spec(\tk).$

\bcor\label{Thm.Gamma-AP.with.parameters}
Let $\Ga$ be a directed rooted tree.
\bee
\item
The property  APP.$(\Ga,t)$  holds  in the Nash case.
\item
The property   SAP.$(\Ga,t)$   holds in the Nash and formal cases, assuming the base field is $\aleph_0$-complete. 
\eee
\ecor
\bpr
Adjoin to the quiver $\Ga$ the additional vertex $ \{Spec(\k)\},$ and one directed edge, $\Ga\ni o\to \{Spec(\k)\}.$
 The ``extended" quiver  $\Ga'=\Ga\cup \{Spec(\k)\} $ is still a directed rooted graph without cycles.

  Repeat the proof of Theorem \ref{Thm.Gamma-AP.for.k<x>}, but in the specializations of Step 2 do not specialize the $t$-parameters.
\epr

\bex\label{Ex.normal.form.unfolding}
The simplest quivers $\Phi_X\circlearrowright X\to Y$ and $\Phi_X\circlearrowright X\to Y\circlearrowleft \Phi_Y$  are particularly important for the study of families/unfoldings of maps, \cite{Kerner.Unfoldings}. Suppose the target is smooth, $Y=(\k^m,o),$ then we identify $\Maps\cong \cm_X\cdot \RmX.$ Each map $f:X\to Y$ is an unfolding of its central fibre, $f_o:X_o\to Y_o.$
 Fix a finite set of vectors, $v_\bullet\in \RmX.$

  We want to bring an unfolding $f_t$ to the ``normal form", $f_o+\sum c_\bullet(t)\cdot v_\bullet.$
  Corollary \ref{Thm.Gamma-AP.with.parameters} implies: if   $f_t$ is $\hat\cR$ (resp. $\widehat{\cL\cR}$) equivalent to
   $f_o+\sum \hc_\bullet(t)\cdot v_\bullet$ (with formal coefficients $c_\bullet\in \k[\![t]\!]$),
  then $f_t$ is
   $\cR$ (resp. $\cL\cR$) equivalent to   $f_o+\sum c_\bullet(t)\cdot v_\bullet$ (with Nash coefficients $c_\bullet\in \k\bl t\br$).
%\vspace{-0.7cm}
\eex


\begin{thebibliography}{\!\!\!\!\!\!99\!\!\!\!}
\bibitem[Abhyankar-van der Put.70]{Abhyankar-van der Put.70}
 S. S. Abhyankar,  M. van der Put, {\em Homomorphisms of analytic local rings.} J. Reine Angew. Math. 242 (1970), 26–-60.


\bibitem[Alon.Cas-Jim.Haus.Kout.18]{Alonso.Castro.Hauser.Koutschan}
M. E. Alonso, F. J. Castro-Jim\'{e}énez, H. Hauser, C. Koutschan, {\em
Echelons of power series and Gabrielov’s counterexample to nested linear Artin approximation.}
Bull. Lond. Math. Soc. 50, No. 4, 649--662 (2018).

\bibitem[A.G.V.]{AGV}  V. I. Arnold, S. M. Gusein-Zade, A. N. Varchenko, {\em
Singularities of differentiable maps. Volume 1. Classification of critical points, caustics and wave fronts.}
   Reprint of the 1985 edition. Modern Birkh\"auser Classics. Birkh\"auser/Springer, New York, 2012


\bibitem[A.G.L.V.]{AGLV} V. I. Arnold, V.V. Goryunov, O.V. Lyashko, V.A. Vasil'ev, {\em Singularity theory. I.}  Reprint of the original English edition from the series Encyclopaedia of Mathematical Sciences . Springer-Verlag, Berlin,1998.iv+245 pp


\bibitem[Artin.68]{Artin.68} M.Artin, {\em On the solutions of analytic equations}.  Invent. Math.  5,  1968, 277--291.

\bibitem[Artin.69]{Artin.69}M. Artin, {\em  Algebraic approximation of structures over complete local rings}, Publ.Math.IHES, 36, (1969), 23-58.

\bibitem[BDLvdD79]{BDLvdD79} J. Becker, J. Denef, L. Lipshitz, L. van den Dries, {\em Ultraproducts and approximation in local rings
I}, Invent. Math., 51, (1979), 189-203.


\bibitem[B.K.16]{B.K.motor}
G. Belitski, D. Kerner,  {\em Group actions on filtered modules and finite determinacy. Finding large submodules
in the orbit by linearization}, C. R. Math. Acad. Sci. Soc. R. Can. 38 (2016),
no. 4, 113--153.


\bibitem[Bel.Boix.Ker.20]{Belitskii.Boix.Kerner}  G. Belitskii, A. F. Boix, D. Kerner, {\em
Approximation results of Artin-Tougeron-type for general filtrations and for $C^r$-equations.} J. Pure Appl. Algebra 224, No. 12, 106431, 14 p. (2020).



%    \bibitem[Belot.Curmi.Rond.21]{Bel.Cur.Ron.}
%  A.  Belotto da Silva, O. Curmi,  G. Rond, {\em A proof of A. Gabrielov's rank theorem.} J. \'{E}c. polytech. Math. 8 (2021), 1329--1396.

\bibitem[Bil.Par.Ron.17]{Bilski.Parusiński.Rond} M.  Bilski, A. Parusi\'{n}ński, G. Rond, {\em
Local topological algebraicity of analytic function germs.} J. Algebr. Geom. 26, No. 1, 177-197 (2017).

\bibitem[B.G.K.22]{BGK.20} A.-F. Boix, G.-M. Greuel, D. Kerner, {\em  Pairs of Lie-type and large orbits of group actions on filtered modules.
 (A characteristic-free approach to finite determinacy.)},  Math. Z. 301 (2022), no. 3, 2415--2463.

%\bibitem[Bourbaki.CA]{Bourbaki.CA} N. Bourbaki, {\em Alg\`{e}bre commutative},  Fasc. XXVIII, Chap. III, \S4, No. 5 (1961).


\bibitem[Bourbaki.Lie]{Bourbaki.Lie} N. Bourbaki, {\em Elements of Mathematics. Lie groups and Lie algebras. Chapters 1--3},
Hermann, Paris; Addison-Wesley Publishing Co., Reading, Mass. 1975, xxviii+450.

%\bibitem[Cas-Jim.Pop.Ron.19]{Cas-Jim.Pop.Ron}
% F.-J. Castro-Jim\'{e}nez, D. Popescu, G. Rond, {\em
%Linear nested Artin approximation theorem for algebraic power series.} Manuscr. Math. 158, No. 1-2, 55-73 (2019).



\bibitem[Denef-Lipshitz.80]{Denef-Lipshitz} J. Denef, L. Lipshitz, {\em Ultraproducts and approximation in local rings. II.}
Math. Ann. 253 (1980), no. 1, 1--28.

%\bibitem[Eisenbud]{Eisenbud-book} D. Eisenbud, {\em Commutative algebra. With a view toward algebraic geometry.}
% Graduate Texts in Mathematics, 150. Springer-Verlag, New York, 1995.


\bibitem[Fich.Quar.Shiot.16]{Fichou.Quarez.Shiota} G. Fichou, R. Quarez, M Shiota, {\em Artin approximation compatible with a change of variables.}
 Canad. Math. Bull. 59 (2016), no. 4, 760–-76.



\bibitem[Gabrielov.71]{Gabrielov.71} A. M.    Gabrielov, {\em Formal relations among analytic functions.}
Izv. Akad. Nauk SSSR Ser. Mat. 37 (1973), 1056--1090.
 Translation from Funkts. Anal. Prilozh. 5, No. 4, 64-65 (1971).


\bibitem[Gabrielov.73]{Gabrielov.73} A. M. Gabrielov, {\em Formale Relationen zwischen analytischen Funktionen.}
Izv. Akad. Nauk SSSR, Ser. Mat. 37, 1056--1090, (1973).




\bibitem[Gr.Lo.Sh]{Gr.Lo.Sh} G.-M. Greuel, C. Lossen, E. Shustin, {\em Introduction to singularities and deformations.}
Second Edition, Springer Monographs in Mathematics. Springer, Berlin, 2024-2025.


%\bibitem[Grothendieck.61]{Grothendieck.61} A. Grothendieck, {\em
%Techniques de construction en g\'eéom\'eétrie analytique. VI. \'EÉtude locale des morphismes: germes d'espaces analytiques, platitude, morphismes %simples}, S\'eéminaire Henri Cartan, Volume 13 (1960-1961) p. 1-13.


%\bibitem[Huneke-Swanson]{Huneke-Swanson}
%C. Huneke, I. Swanson, {\em Integral closure of ideals, rings, and modules},
%London Math. Soc. Lecture Note Ser., 336
%Cambridge University Press, Cambridge, 2006, xiv+431 pp.



\bibitem[Izumi.07]{Izumi.07} S. Izumi, {\em Fundamental properties of germs of analytic mappings of analytic sets and related topics.}
  `Real and complex singularities'. Hackensack, NJ: World Scientific, 109-123 (2007).



%\bibitem[Joi\c{t}ţa-Tibar.21]{Joiţa-Tibar} C. Joi\c{t}ţa, M. Tib\u{a}ăr, {\em The local image problem for complex analytic maps},
%Ark. Mat. 59 (2021), no. 2, 345-–358.


%\bibitem[K\"aällstr\"oöm.13]{Kallstrom}   R. K\"aällstr\"oöm, {\em Purity of branch and critical locus,}
%J. Algebra 379 (2013), 156-–178.



\bibitem[Kerner.24]{Kerner.Group.Orbits} D. Kerner, {\em  Orbits of the left-right equivalence of maps in arbitrary characteristic},  European Journal of Mathematics 10, 64 (2024).


\bibitem[Kerner.22]{Kerner.Unfoldings} D. Kerner, {\em
Unfoldings of maps, the first results on stable maps, and results of Mather-Yau/Gaffney-Hauser type in arbitrary characteristic},
arXiv:2209.05071.

\bibitem[Kerner.25]{Kerner.LRAP} D. Kerner, {\em Results on left-right Artin approximation for algebraic morphisms and for analytic morphisms of weakly-finite singularity type}, J. Lond. Math. Soc. (2) 111 (2025), no. 1, Paper No. e70053.


\bibitem[Kerner.25.prep]{Kerner.Change.of.Base.Field} D. Kerner, {\em Equivalence of mapping-germs over $\k$ vs that over $\K$}, preprint.


\bibitem[Kirillov]{Kirillov} A. Kirillov, {\em Quiver representations and quiver varieties.} 
Graduate Studies in Mathematics 174. Providence, RI: American Mathematical Society (AMS)  xii, 295 p. (2016).

\bibitem[Kosar.Popescu]{Kosar.Popescu} Z. Kosar, D. Popescu, {\em 
Nested Artin strong approximation property.} 
J. Pure Appl. Algebra 222, No. 4, 818-827 (2018).

%\bibitem[K.P.P.R.M.78]{K.P.P.R.M.} H. Kurke, G. Pfister, D. Popescu, M. Roczen, T. Mostowski, {\em Die Approximationseigenschaft lokaler Ringe,} %Lectures Notes in Math., 634, Springer-Verlag, Berlin-New York, (1978).


\bibitem[Matsumura]{Matsumura} H. Matsumura, {\em Commutative ring theory,}
Cambridge Studies in Advanced Mathematics, 8. Cambridge etc.: Cambridge University Press. XIII, 320 p.

%\bibitem[Milnor]{Milnor} J. Milnor, {\em  Dynamics in one complex variable,}
%Friedr. Vieweg \& Sohn, Braunschweig, 1999,  257 pp.

\bibitem[Milnor]{Milnor} J. Milnor, {\em  Dynamics in one complex variable,}
Friedr. Vieweg \& Sohn, Braunschweig, 1999,  257 pp.


\bibitem[Mond-Montaldi.91]{Mond-Montaldi} D. Mond, J. Montaldi, {\em Deformations of maps on complete intersections, Damon's $\cK_V$-equivalence and bifurcations.} Brasselet, Jean-Paul (ed.),  `Singularities in geometry and topology'.   Cambridge: Cambridge University Press. Lond. Math. Soc. Lect. Note Ser. 201, 263-284 (1994).


\bibitem[Mond. Nu\~{n}.-Bal.]{Mond-Nuno} D. Mond, J.-J. Nu\~{n}o-Ballesteros, {\em Singularities of Mappings.
The Local Behaviour of Smooth and Complex Analytic Mappings}, Grundlehren Math. Wiss., 357,
Springer, Cham, 2020, xv+567 pp.

%\bibitem[Moret-Bailly.21]{Moret-Bailly}     L. Moret-Bailly, {\em     A Henselian preparation theorem},   Israel Journal of Mathematics, 257, %519-531 (2023).


\bibitem[Ogoma.94]{Ogoma.94} T. Ogoma, {\em General N\'eéron desingularization based on the idea of Popescu,}
J. Algebra 167 (1994), no. 1, 57–-84.

%\bibitem[Osgood.1916]{Osgood} W. Osgood, {\em On functions of several complex variables}, Trans. Amer. Math. Soc. 17 (1916), 1-8.

\bibitem[Pfister-Popescu.75]{Pfister-Popescu} G.  Pfister, D. Popescu, {\em Die strenge Approximationseigenschaft lokaler Ringe.}
Invent. Math. 30 (1975), no. 2, 145--174.



\bibitem[P\l oski.74]{Ploski1974} A. P\l oski, {\em Note on a theorem of M. Artin,}  Bull. Acad. Polon. Sci. S\'er. Sci. Math. Astronom. Phys., 22, (1974), 1107--1109.


\bibitem[Popescu.Rond]{Popescu.Rond} D. Popescu, G. Rond, {\em  Remarks on Artin approximation with constraints.}
Osaka J. Math. 56, No. 3, 431-440 (2019)

\bibitem[Popescu.85]{Popescu.85} D. Popescu, {\em General N\'{e}éron desingularization.} Nagoya Math. J. 100, 97--126 (1985).



\bibitem[Popescu.86]{Popescu.86} D. Popescu, {\em General N\'{e}éron desingularization and approximation.} Nagoya Math. J.  104, 85--115 (1986).





\bibitem[Popescu.00]{Popescu} D. Popescu,
{\em  Commutative Rings and Algebras. Artin approximation},   Handbook of Algebra, Vol 2, pages 321-356, 2000.


\bibitem[Popescu-Rond]{Popescu-Rond} D. Popescu, G. Rond, {\em  Remarks on Artin approximation with constraints.}
Osaka J. Math. 56, No. 3, 431-440 (2019).


\bibitem[Rond.18]{Rond}   G. Rond, {\em Artin approximation},  J. Singul. 17 (2018), 108--192.


\bibitem[Shiota.98]{Shiota.1998} M. Shiota, {\em Relation between equivalence relations of maps and functions.}
  Real analytic and algebraic singularities, 114--144, Pitman Res. Notes Math. Ser., 381, Longman, Harlow, 1998.

\bibitem[Shiota.10]{Shiota.2010} M. Shiota, {\em Analytic and Nash equivalence relations of Nash maps.} Bull. Lond. Math. Soc. 42 (2010), no. 6, 1055--1064.



\bibitem[Spivakovsky.99]{Spivakovsky.99} M. Spivakovsky, {\em A new proof of D. Popescu's theorem on smoothing of ring homomorphisms},   J. Amer. Math. Soc., 12  (1999),  no. 2, 381-444.


\bibitem[Swan.98]{Swan.98} R. G. Swan, {\em N\'eéron-Popescu desingularization} Lect. Algebra Geom., 2
International Press, Cambridge, MA, 1998, 135–-192.
\end{thebibliography}
\end{document}